\documentclass[12pt,leqno]{article}

\usepackage[T1]{fontenc}
\usepackage[utf8]{inputenc}

\usepackage{amsmath}
\usepackage{amssymb}
\usepackage[mathscr]{eucal}
\usepackage{graphicx}
\usepackage{color}

\usepackage{caption}
\usepackage{subcaption}
\usepackage{float}

\newtheorem{Theorem}{\sc Theorem}

\newtheorem{Lemma}[Theorem]{\sc Lemma}

\newtheorem{Example}[Theorem]{\sc Example}

\newtheorem{Experiment}{\sc Experiment}

\newcommand{\R}{{\if mm {\rm I}\mkern -3mu{\rm R}\else \leavevmode
		\hbox{I}\kern -.17em\hbox{R} \fi}}

\newcommand{\cF}{\mbox{{${\cal F}$}}}

\newcommand{\cQ}{\mbox{{${\cal Q}$}}}

\newcommand{\cP}{\mbox{{${\cal P}$}}}

\newcounter{theorem}

\def\sqr#1#2{{
		\vcenter{
			\vbox{\hrule height.#2pt
				\hbox{\vrule width.#2pt height#1pt \kern#1pt
					\vrule width.#2pt
				}
				\hrule height.#2pt
			}
		}
}}
\def\square{\mathchoice\sqr84\sqr84\sqr{2.1}3\sqr{1.5}3}

\def\real{\mathbb{R}}

\def\lista#1
{{ \itemindent 0.0cm \labelsep .2cm \leftmargin 0.8cm \rightmargin
		0.0cm \labelwidth 0.6cm \topsep 0.0mm
		\parsep 0.0mm
		\itemsep 0.0mm
		\begin{list}{}
			{ \setlength{\leftmargin}{.8cm} \setlength{\rightmargin}{0.0cm}
				\setlength{\parsep}{0.0mm} \setlength{\topsep}{.0mm}
				\setlength{\parskip}{.0cm} \setlength{\itemsep}{.0cm} }
			{#1}\end{list}} }

\begin{document}
	
	\title{\bf Modelling, Analysis and Numerical Simulation of a  Spring-Rods System with Unilateral Constraints}

	\vspace{5mm}
	

	\vspace{22mm}  
	{\author{Anna Ochal$^{1}$,\ Wiktor Prz\k{a}dka$^{1}$,\ Mircea Sofonea$^{2}$\footnote{Corresponding author, E-mail : sofonea@univ-perp.fr}\\ and\\ Domingo A. Tarzia$^{3,4}$\\[4mm]	
			{\it \small $1$ Chair of Optimization and Control, }
			{\it \small
				Jagiellonian University in Krakow}
			\\{\it\small Ul. Lojasiewicza 6, 30348 Krakow, Poland}		\\[4mm]	
			{\it \small $^2$ Laboratoire de Math\'ematiques et Physique}
			{\it \small
				University of Perpignan Via Domitia}
			\\{\it\small 52 Avenue Paul Alduy, 66860 Perpignan, France}		\\[5mm]		
			{\it\small $^{3}$  Departamento de Matem\'atica, FCE}\\ {\it \small Universidad Austral}
			{\it \small Paraguay 1950, S2000FZF Rosario, Argentina}\\[4mm]
		{\it\small $^{4}$ CONICET, Argentina}}


	\date{}
	\maketitle
	\thispagestyle{empty}
	
	\vskip 0mm
	
	\noindent {\bf Abstract.} \ 
	In this paper we  consider a mathematical model which describes the equilibrium of two elastic rods attached to a nonlinear spring. We derive the variational formulation of the model which is in the form of an elliptic quasivariational inequality for the displacement field. We prove the unique weak solvability of the problem, then we state and prove some convergence results, for which we provide the corresponding mechanical interpretation. Next, we turn to the numerical approximation of the problem based on a finite element scheme. We use a relaxation method to solve the discrete problems that we implement on the computer. Using this method, we provide numerical simulations which validate our convergence results.

	\vskip 2mm 
	
	
	\vskip 2mm
	

	
	
	
	
	

	\vskip 15mm
	
	\section{Introduction}\label{s1}
	\setcounter{equation}0

	\vskip2mm\noindent
	The interest in mathematical problems that describe the contact of thin structures  like spring, rods and beams is two folds. First, such kind of problems arise in
	real-word setting like automotive industry and, more specifically, in  motors, engines and transmissions; they represent an important kind of problems in Mechanics of Structures, too. 
	Second, the study of these problems represents a first step in the modeling and analysis of  more complicate problems which describe the contact of deformable bodies in the three-dimensional setting.  Even stated in a simplified framework, the study of mechanical systems involving spring, rods and beams gives rise to interesting mathematical models, which are interesting in their own. The literature in the field includes \cite{AS, AFS, ABS, CRSS, G, HSU, SB}, for instance. There, various models of spring-rod systems have been considered, constructed by using different contact conditions and interface laws, and unique solvability results have been proved. 
	
	This current paper is dedicated to the modeling, analysis and numerical simulation of a new spring-rods system. The novelty arises in the fact that here we consider a mechanical system made by two nonlinear elastic rods connected at their ends with an elastic spring which could be completely squeezed and, therefore, the rods could arrive in contact; the contact is unilateral and it is described by using a Signorini-type condition. Considering this physical setting leads to a non standard mathematical model, stated in terms of an elliptic variational inequality formulated in a product Hilbert space. The unknown of the model is the couple in which the components are
	the displacement fields in the two rods.
	Existence and uniqueness results on variational inequalities can be found in \cite{BC, B, Br, Kind-St,L}, for instance. They are obtained by using
	arguments of monotonicity and convexity, including properties of the subdifferential of a convex function.  The numerical analysis of various classes of variational inequalities was treated in \cite{Gl,GLT,HHNL}, among others. The theory of variational inequalities finds various applications
	in Mechanics and Physics and, in particular, in Contact Mechanics, as illustrated in \cite {DL,EJK,HR, HS,HHN,KO, P,SM1,SM2}, for instance.

	The rest of the paper is structured as follows.  In Section \ref{s1n} we present the mathematical preliminaries we need in the analysis of our contact problem. These preliminaries concern existence, uniqueness and convergence results obtained in \cite{SM2,ST2} in the study of elliptic quasivatiational inequalities. In Section
    \ref{s2} we present the physical setting of the spring-rods system we consider, then we list the mechanical assumptions and state the corresponding mathematical model. In Section
	\ref{s3} we derive a variational formulation of the model, obtained by
	using a standard procedure based on integration by parts.
	Then we provide the existence of a unique weak solution to the problem. In Section \ref{s4} we study the behavior of the solution when the stiffness of the spring converges to infinite. Here, we state and prove three convergence results and provide their mechanical interpretation. In Section
	\ref{s5} we introduce a finite-dimensional discrete scheme to approach our model, 
	describe the algorithm we use to solve the corresponding discrete problems, then we present numerical simulations and provide the corresponding mechanical interpretations. The simulations we provide here represent an evidence of our theoretical convergence results.

\section{Preliminaries}\label{s1n}
	\setcounter{equation}0

In this section, we recall some results in the study of a general class of elliptic quasivariational inequalities. The functional framework is the following: $X$ is a real Hilbert space endowed with the inner product $(\cdot,\cdot)_X$ and norm $\|\cdot\|_X$, $K\subset X$, $A\colon X\to X$, $j\colon X\times X\to\R$ and $f\in X$. Then, the inequality problem we consider is the following.

\medskip\medskip\noindent{\bf Problem}  ${\cal Q}$. {\it Find $u$ such that}
\begin{equation}\label{1a}u\in K ,\qquad(Au,v-u)_X+j(u,v)-j(u,u) \ge(f,v-u)_X \qquad\forall\,v\in K.
\end{equation}
		
	
	\medskip
	In the study of Problem $\cQ$ we consider the following assumptions.
	\begin{eqnarray}
		&&\label{Ka}
		\quad\ \ K \ \mbox{\rm is a nonempty, closed, convex subset of} \ X.
		\\[2.5mm]
		&&\label{Aa}
		\left\{ \begin{array}{l}
			A\ {\rm is\ a\ strongly\ monotone\ Lipschitz\ continuous\ operator,\, i.e.,}\\[0mm] 
			{\rm there\ exist}\ m>0\ {\rm and}\ M>0\ {\rm such\ that}\quad\ \ \\[2mm]
			\mbox{\rm (a) }\ (Au - Av,u -v)_{X} \geq m \|u -v\|^{2}_{X}
			\quad \forall\,u,\,v\in X,\\[2mm]
			\mbox{\rm (b) }\ \|Au-Av\|_X\le M\, {\|u-v\|_X}\quad\forall\,u,\,v\in
			X.
		\end{array}\right.\\[2.5mm]
		&&\label{ja}\left\{\begin{array}{l}  
			\mbox{(a) } {\rm For\ all}\ \eta\in X,\ j(\eta,\cdot)\colon X\to \mathbb
			R \mbox{ is
				convex  and lower semicontinuous.}
               \\[2mm]
			\mbox{(b) There exists  }\alpha\ge 0 \mbox{ such that }\\
			\quad\quad
			j(\eta,\widetilde{v})-j(\eta,v)+
			 j(\widetilde{\eta},v)-j(\widetilde{\eta},\widetilde{v})\\
			\qquad\quad\leq \alpha\,\|\eta-\widetilde{\eta}\|_X\|v-\widetilde{v}\|_X\quad
			\forall \,\eta,\,\widetilde{\eta},\,v,\,\widetilde{v}\in X. \end{array}\right.\\[2.5mm]
		&&\label{small}
		\qquad\quad\ m>\alpha.
		\\ [3mm]
		&&\qquad\quad\label{fa}f\in X.
	\end{eqnarray}

	\medskip
	A proof of the following existence and uniqueness result can be found in \cite{SM2}.
	
	\begin{Theorem}\label{t0a} 
		Assume $(\ref{Ka})$--$(\ref{fa})$.
		Then Problem $\cQ$ has a unique solution.
	\end{Theorem}
	
Consider now a sequence $\{\lambda_n\}\subset\real$ and an operator $G\colon X\to X$. For each $n\in\mathbb{N}$
denote by $\cQ_n$  the following version of Problem $\cQ$.

\medskip\noindent{\bf Problem}  ${\cal Q}_n$. {\it Find $u_n$ such that}
\begin{eqnarray}\label{1n}
	&&u_n\in K,\quad (A u_n, v - u_n)_X +
	\frac{1}{\lambda_n} (Gu_n, v - u_n)_X +
	j(u_n,v) - j(u_n,u_n)\label{2n}\\[2mm]
	&&\qquad \qquad\qquad
	\ge (f_n, v - u_n)_X
	\qquad\forall\, v \in K.\nonumber
	\nonumber
\end{eqnarray}

Moreover, let $K^*\subset X$, and consider the following inequality problem in which the constraints are governed by the set $K^*$.

\medskip\medskip\noindent{\bf Problem}  ${\cal Q}^*$. {\it Find $u^*$ such that}
\begin{equation}\label{1b}u^*\in K^* ,\qquad(Au^*,v-u^*)_X+j(u^*,v)-j(u^*,u^*) \ge(f,v-u^*)_X \qquad\forall\,v\in K^*.
\end{equation}

In the study of these problems, we consider the following assumptions.
\begin{eqnarray}
	&&\label{za1}
	K^* \ \mbox{\rm is a nonempty, closed, convex subset of $X$.}\\ [2mm]
	&&\label{za2} K^*\subset K. \\[2mm]
	&&\label{za3}G:X\to X \ \ \mbox{is a  monotone  Lipschitz continuous operator}.\\ [2mm]
	&&\label{za4}\left\{\begin{array}{ll}
		\mbox{\rm (a)}
		\quad(Gu,v-u)_X\le 0\qquad\forall\, u\in K,\ v\in K^*.\\ [3mm]
		\mbox{\rm (b)} \quad u\in K,\quad (Gu,v-u)_X=0\quad\forall\,v\in K^*\ \ \Longrightarrow\ \  u\in K^*.
	\end{array}\right.\\ [2mm]
	&&\label{za5}\lambda_n>0\qquad\forall\,n\in\mathbb{N}. \\ [2mm]
	&&\label{za6}\lambda_n\to 0\quad {\rm as}\quad n\to\infty.
	\end{eqnarray}

\medskip
We have the following  existence,  uniqueness
and convergence result. 

\begin{Theorem}\label{t1a}
	Assume  $(\ref{Ka})$--$(\ref{fa})$ and $(\ref{za1})$--$(\ref{za6})$.
	Then, Problem $\cQ^*$ has a unique solution $u^*\in K^*$ and, for each $n\in\mathbb{N}$, there exists a unique solution $u_n\in K$
	to Problem~${\cal Q}_n$. Moreover,
	$u_n\to u^*$ in $X$, as $n\to \infty$.
\end{Theorem}

\medskip
A proof of the theorem can be found in \cite{ST2}.
Note that  the existence and uniqueness part in this theorem is a direct consequence of Theorem \ref{t0a}. The  convergence part was obtained in several steps, by using arguments of compactness, pseudomonotonicity and lower semicontinuity.

A brief comparison between inequalities (\ref{1b}) and (\ref{2n}) shows that   (\ref{2n}) is obtained from (\ref{1b}) by replacing the set $K^*$ with a larger set  $K$ and the operator $A$ with the operator $A+\frac{1}{\lambda_n}\,G$, in which $\lambda_n$ is a penalty parameter.  For this reason we refer to (\ref{2n})  as a penalty problem of  (\ref{1b}). Theorem \ref{t1a} establishes the link between the solutions of these problems.  Roughly speaking, it shows that, in the limit when $n\to\infty$, a partial relaxation of the set of constraints can be compensated by a convenient perturbation of the nonlinear operator which governs Problem $\cQ^*$.
Finally, note  that
Problem~${\cal Q}_n$ represents a penalty problem of $\cQ^*$.
Penalty methods have been widely used in the literature as an approximation tool to treat constraints in variational and hemivariational inequalities, as explained in \cite{G,KO,SM2,SMBOOK} and the references therein.

 \medskip We end this section with an example of operator $G$ which satisfies conditions (\ref{za3}) and  (\ref{za4}).

 \begin{Example} Assume that $\eqref{Ka}$,
 	$\eqref{za1}$ and $\eqref{za2}$  hold and denote by ${P}_{K^*} \colon X \to K^*$ and  ${P}_{K}\colon X \to K$  the projection operators on the sets $K^*$ and ${K}$, respectively. Then, using the properties of the projection operators it is easy to see that the operator $G:X\to X$ given by  $G=2I_X-P_{K^*}-P_{K}$ satisfies conditions $(\ref{za3})$ and  $(\ref{za4})$.
 \end{Example}

	\section{The model}\label{s2}
	\setcounter{equation}0
	
	In this section we introduce the physical setting, then we construct the corresponding mathematical model which describes the equilibrium of the spring-rods system. 
	
The physical setting is as follows. 	We consider two elastic rods which, in their reference configuration, occupy  the intervals $[a,-l]$  and   $[l,b]$  on the $Ox$ axis, respectively. Here $a<0$, $b>0$ and $l>0$ are given constants such that $a<-l$ and $l<b$. Therefore, before the deformation process, the length of the first rod is $L_1=-l-a>0$ and   the length of the second rod is $L_2=b-l>0$. The rods are fixed  at their  ends $x=a$  and $x=b$ and their extremities $x=-l$ and $x=l$ are attached to a nonlinear spring.  The natural length of the spring is $2l$ and, in the reference configuration,  no forces are acting on the rod. This situation corresponds to Figure~\ref{fig}a).

	Assume now that the rods are submitted to the action of body forces of {line} density $f_1$ and $f_2$ which act along the $Ox$ axis. As a result, the mechanical system evolves to an equilibrium configuration in which the rods and the spring are deformed. In this configuration
	the spring could be either in extension (as depicted in Figure~\ref{fig}b))
	or in compression  (as depicted in Figure~\ref{fig}c)). We assume that the spring has an elastic behaviour and could be completely compressed. When this situation arises, the current length of the spring vanishes and the two rods arrive in contact. The contact of the two ends of the rods is without penetration.
	
	Denote  by $u_1$, $u_2$ the displacement field in the two rods and by $\sigma_1$, $\sigma_2$ the corresponding stress fields, respectively. 
	Then, the problem of finding the equilibrium of the mechanical system  in the physical setting described above can be formulated as follows.

	\vspace{00mm}
	
	\begin{center}
		\begin{figure}[htbp] 
			\hspace{10mm}\includegraphics[width=5in]{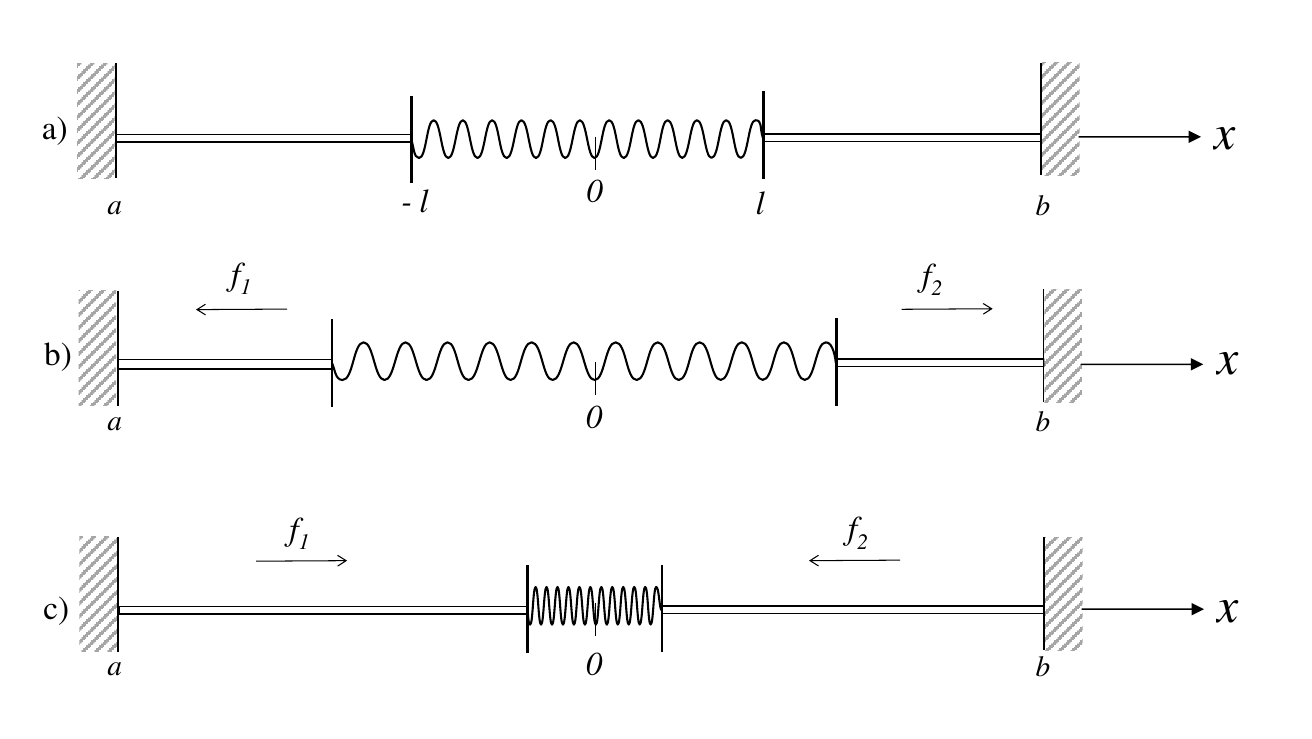}
			\vspace{00mm}
			\caption{The spring-rods system with unilateral constraints:}
			
			\centerline {	a) reference configuration;\quad b) spring in extension;\quad c) spring in compression.}\label{fig}
		\end{figure}
	\end{center}

	\medskip
	\noindent
	{\bf Problem} ${\cal P}$. {\it Find the displacement fields  $u_1\colon [a,-l]\to\mathbb{R}$,
		$u_2\colon [l,b]\to\mathbb{R}$  and the stress fields $\sigma_1\colon [a,-l]\to\mathbb{R}$, $\sigma_2\colon [l,b]\to\mathbb{R}$ such that} 
	     \begin{eqnarray}
		&&\quad  \sigma_1(x) = {\cal F}_1\Big(x,\frac{du_1}{dx}(x)\Big) \quad \mbox{for}\ x\in\ (a,-l)
		,\label{2.1}\\[3mm]
		&&\quad  \sigma_2(x) = {\cal F}_2\Big(x,\frac{du_2}{dx}(x)\Big) \quad \mbox{for}\ x\in\ (l,b)
		,\label{2.1n}\\[3mm]
		&&\quad  \frac{d\sigma_1}{dx}(x) + f_1(x) = 0\quad \ \ \ \ \  \mbox{for}\  x\in(a,-l),
		\label{2.2}\\[3mm]
		&&\quad  \frac{d\sigma_2}{dx}(x) + f_2(x) = 0\quad \ \ \ \ \ \mbox{for}\  x\in(l,b),
		\label{2.2n}\\[3mm]
		&&\label{2.3}\quad \ u_1(a) = 0, \\[3mm]
		&&\label{2.3n}\quad \ u_2(b) = 0, \\[3mm]
		&&\left\{ 
		\begin{array}{ll}
			u_1(-l)-u_2(l)\le2l,\\[2mm]
			\sigma_1(-l)=\sigma_2(l)=-p\big(2l+u_2(l)-u_1(-l)\big)
            &{\rm if}\ u_1(-l)-u_2(l)<2l,\\[2mm]
			\sigma_1(-l)=\sigma_2(l)\le-p\big(2l+u_2(l)-u_1(-l)\big){=-p(0)}
            &{\rm if}\ u_1(-l)-u_2(l)=2l.
			\label{2.4}\end{array}\right.
	\end{eqnarray}
	
	\medskip
	
	We now provide a short description of the equations and boundary condition in Problem  ${\cal P}$. 
	
	First, equations (\ref{2.1}) and (\ref{2.1n}) represent the elastic constitutive laws of the rods in which  the derivatives $\frac{du_i}{dx}$  represent the linear strain field and ${\cal F}_i:\mathbb{R}\to\mathbb{R}$ are the elasticity operators, for $i=1,2$. 
	Equations (\ref{2.2}) and (\ref{2.2n}) are the balance equations which describe the equilibrium of the rods,  and conditions (\ref{2.3}), (\ref{2.3n}) represent the displacement conditions at the outer ends. We use these conditions here since the rods are assumed to be fixed  at the ends $x=a$ and $x=b$, respectively. 
	
	Finally, condition (\ref{2.4}) represents the contact boundary condition in which   $p$ is a given real-valued function which models the behaviour of the spring. It could be nonlinear and it will be described below. 
	We assume in what follows that the  spring is not damageable and, therefore,  the function $p$ in (\ref{2.4}) has a unique argument.  Considering the case of a damageable spring could be possible.

 Nevertheless, as explained in \cite{CRSS},  in this case the problem becomes evolutionary and the function $p$ depends also on a new time-dependent variable, the damage function.
	Our interest lies in condition (\ref{2.4}) and, therefore, we describe it below, line by line.

	To this end, recall that  the exterior unit normal
	at the end $x=-l$ of the first rod is $\nu_1=1$ and, therefore $\sigma_1(-l)$ represents the stress vector at the point $x=-l$ of the first rod. In contrast, the exterior unit normal
	at the extremity $x=l$ of the second rod is $\nu_2=-1$ and, therefore $\sigma_2(l)$ represents the opposite of the stress vector at the extremity $x=l$ of the second rod. We need this remark in the arguments we provide below, in order to 
	describe the sense of these vectors, which could be towards the origin or towards the rods.

	Next, note that the quantity $\theta=2l+u_2(l)-u_1(-l)$ represents the current length of the spring and, therefore, condition 	$u_1(-l)-u_2(l)\le2l$  shows that this length is always non-negative. Indeed, if $\theta<0$ then the rods would penetrate themselves which is not physically accepted. To conclude, the first inequality in (\ref{2.4}) represents the nonpenetrability condition.
	
	Assume now that $\theta=2l+u_2(l)-u_1(l)>2l$, i.e., the spring is in extension.
	Then, since the spring is assumed to be elastic, it
	exerts a traction on the ends of two rods which depends on its length $\theta$, and which is
	towards of the origin $O$ of the system.
	Thus, we have $\sigma_1(-l)=\sigma_2(l)=-p(\theta)$
	and  $p$ has to be such that $p(r)<0$ for $r>2l$, in order to guarantee that $\sigma_1(-l)>0$ and $-\sigma_2(l)<0$, i.e., the spring pulls the rods.
	When $0<\theta=2l+u_2(l)-u_1(l)<2l$ then
	the spring is in compression. We still have
	$\sigma_1(-l)=\sigma_2(l)=-p(\theta)$ and we impose that $p(r)>0$ for $r<2l$ 
	in order to guarantee that $\sigma_1(-l)<0$ and $-\sigma_2(l)>0$. 
   These inequalities show that in this case the reaction of the spring is towards the rods, i.e., the spring pushes the rods.
	These arguments justify the second line in the contact condition (\ref{2.4}). 
	
	Assume now that $\theta=2l+u_2(l)-u_1(l)=0$. Then,  the spring is fully compressed and two extremities of the rods arrive in contact. In this case,  
	besides the {pressure} $p(0)$ exerted by the compressed spring, at the ends $x=-l$ and $x=l$ of the rods we have an additional pressure which models the reaction of each rod on the other one. These reactions are  towards the rods and prevent the material penetrability. More precisely, the additional {pressure} which acts on the first rod is  negative  and it represents the reaction of the second rod. In contrast, the additional {pressure} which acts on the second rod is positive  and it represents the reaction of the first rod. This justifies the third line in (\ref{2.4})
	since, recall, $\sigma_1(-l)$ represents the stress vector at the point $x=-l$ of the first rod while $\sigma_2(l)$ represents the opposite of the stress vector at the extremity $x=l$ of the second rod.
	
	Note that Problem $\cP$ introduced above is nonlinear and involves unilateral constraints. For this reason, its analysis will be done through its variational formulation  that we present in the next section.

	\section{Unique weak solvability}\label{s3}
	\setcounter{equation}0

	\medskip
	In the study of Problem ${\cal P}$ we use the standard notation for Lebesgue and Sobolev spaces. In addition, for the displacement fields $u_1$ and $u_2$ we need the spaces 
	\begin{eqnarray*}
		&&V_1 = \lbrace\, v_1 \in H^1(a,-l): v_1(a)=0\, \rbrace,\\ 
		&&V_2 = \lbrace\, v_2 \in H^1(l,b): v_2(b)=0\, \rbrace. 
	\end{eqnarray*}
	It is well known that
	the spaces  $V_1$ and $V_2$
	are real Hilbert spaces with the canonical inner products 
	\begin{equation*}
		(u_1, v_1)_{V_1}= \int_a^{-l}\frac{du_1}{dx}\,\frac{dv_1}{dx}\,dx,\qquad 
		(u_2, v_2)_{V_2}= \int_l^{b}\frac{du_2}{dx}\,\frac{dv_2}{dx}\,dx 
	\end{equation*}
	for all $u_1,\, v_1\in V_1$ and $u_2,\, v_2\in V_2$, respectively.
	The associated norms on these spaces will be denoted by  $\|\cdot\|_{V_1}$ and $\|\cdot\|_{V_2}$.
	Moreover,  using the identities
	\[v_1(-l)=\int_a^{-l}\frac{dv_1}{dx}\,dx,\quad v_2(l)=-\int_l^{b}\frac{dv_2}{dx}\,dx\]
	and the Cauchy-Schwarz inequality, it follows that
	\begin{equation}\label{inj}
		|v_1(-l)| \leq\sqrt{L_1}\,\|v_1\|_{V_1}, \ \quad |v_2(l)| \leq\sqrt{L_2}\,\|v_2\|_{V_2}  \
	\end{equation}
	for all \, $v_1 \in V_1,\ v_2\in V_2$. Recall that, here and below, $L_1=-l-a>0$ and $L_2=b-l>0$. Therefore, with notation $L=\max\,\{L_1,L_2\}$ inequalities \eqref{inj}  imply that
	\begin{equation}\label{injn}
		|v_1(-l)| \leq\sqrt{L}\,\|v_1\|_{V_1}, \ \quad |v_2(l)| \leq\sqrt{L}\,\|v_2\|_{V_2}  \
	\end{equation}
	for all \, $v_1 \in V_1,\ v_2\in V_2$.

	Let $V$ denote the product of the spaces $V_1$ and $V_2$, that is
	$V=V_1\times V_2$. Then,
	$V$~is a real Hilbert space with the canonical inner product 
	\begin{equation*}
		(u, v)_{V}= \int_a^{-l}\frac{du_1}{dx}\,\frac{dv_1}{dx}\,dx+
		\int_l^{b}\frac{du_2}{dx}\,\frac{dv_2}{dx}\,dx \qquad \forall \ u=(u_1,u_2),\ v=(v_1,v_2) \in V 
	\end{equation*}
	and the associated norm $\|\cdot\|_{V}$.  
	
	\medskip
	We now turn to the variational formulation of Problem ${\cal P}$ and, to this end, we assume that the elasticity operators $\cF_1$ and $\cF_2$ satisfy the following conditions.
	\begin{equation}
		\left\{\begin{array}{ll} 
			{\rm (a)\ }
			\cF_1 \colon (a,-l)\times\mathbb{R}\to\mathbb{R}.\\[2mm]
			{\rm (b)\ There\ exists}\ m_1>0 {\rm\ such\ that\ }\\
			{}\ \qquad  (\cF_1(x,r)-\cF_1(x,s))(r-s)
			\ge m_1\,|r-s|^2\\
			\ \qquad\qquad\forall\,
			r,\,s\in \mathbb{R},\ {\rm a.e.}\ x\in(a,-l).\qquad\qquad\\[2mm]
			{\rm (c)\ There\ exists}\, M_1>0 {\rm\ such\ that\ }\\
			{}\ \qquad  |\cF_1(x,r)-\cF_1(x,s)|
			\le M_1\,|r-s|\\
			\ \qquad\qquad\forall\,
			r,\,s\in \mathbb{R},\ {\rm a.e.}\ x\in(a,-l).\qquad\qquad\\[2mm]
			{\rm (d)}\ \cF_1(x,0)=0\quad {\rm a.e.}\  x\in(a,-l).
		\end{array}\right.
		\label{F1}
	\end{equation}
	\begin{equation}
		\hspace{-2mm}\left\{\begin{array}{ll} 
			{\rm (a)\ }
			\cF_2\colon (l,b)\times\mathbb{R}\to\mathbb{R}.\\[2mm]
			{\rm (b)\ There\ exists}\ m_2>0 {\rm\ such\ that\ }\\
			{}\ \qquad  (\cF_2(x,r)-\cF_2(x,s))(r-s)
			\ge m_2\,|r-s|^2\\
			\ \qquad\qquad\forall\,
			r,\,s\in \mathbb{R},\ {\rm a.e.}\ x\in(l,b).\qquad\qquad\\[2mm]
			{\rm (c)\ There\ exists}\, M_2>0 {\rm\ such\ that\ }\\
			{}\ \qquad  |\cF_2(x,r)-\cF_2(x,s)|
			\le M_2\,|r-s|\\
			\ \qquad\qquad\forall\,
			r,\,s\in \mathbb{R},\ {\rm a.e.}\ x\in(l,b).\qquad\qquad\\[2mm]
			{\rm (d)}\ \cF_2(x,0)=0\quad {\rm a.e.}\  x\in(l,b).
		\end{array}\right.
		\label{F2}
	\end{equation}

	The function $p$ is such that
	\begin{equation}
		\hspace{10mm}	\left\{\begin{array}{ll} {\rm (a)\ }
			p\colon \mathbb{R}\to\mathbb{R}.\\[2mm]
			{\rm (b)\ There\ exists}\, L_p>0 {\rm\ such\ that\ }\\
			{}\ \qquad  |p(r_1)-p(r_2)|
			\le L_p\,|r_1-r_2|\quad\forall\,
			r_1,\,r_2\in \mathbb{R}.\qquad\qquad\\[2mm]
			{\rm (c)\ }p(r)> 0\ \ {\rm if\ }\ r <2l \ \ {\rm and}\ \ p(r)< 0\ \ {\rm if\ }\ r > 2l.
		\end{array}\right.
		\label{pc}
	\end{equation}
	A simple example of nonlinear function $p$ which satisfies condition (\ref{pc}) is given by
	\begin{equation}
		p(r)=
		\left\{\begin{array}{ll}
			-k_1(r-2l)\qquad {\rm if}\quad r < 2l,\\[2mm]
			-k_2(r-2l)\qquad {\rm if}\quad r\ge 2l.\end{array}\right.
		\label{pp}
	\end{equation}	
	Here, $k_1$ and $k_2$ represent the stiffness coefficients of the spring in compression and extension, respectively, and are assumed to be positive. The case $k_1\ne k_2$ arises when the spring is nonlinear and has a different behaviour in extension and compression. The case  $k_1=k_2=k$ corresponds to a linear
	spring of stiffness $k$.

	We also assume that
	\begin{equation}\label{sm}
		m_1+m_2 >2L_pL
	\end{equation}
	and we interpret this condition as a smallness condition for the Lipschitz constant of the function $p$.
	Finally, we assume that
	the densities of  body forces have the regularity
	\begin{equation}\label{fc}
		f_1 \in L^2(a,-l),\qquad f_2 \in L^2(l,b).
	\end{equation}

	Under these assumptions we define the functional $\theta \colon V\to\mathbb{R}$, the set $K$, the operator $A\colon V\rightarrow V$, the function
	$j\colon V\times V\rightarrow \mathbb{R}$ 
	and the element $f$ by equalities
	\begin{eqnarray}
		&&\label{tm}\theta(v)=2l-v_1(-l)+v_2(l)\quad \forall\,\,v=(v_1,v_2)\,\in V,\\ [6mm]
		&&\label{Km}K = \lbrace v=(v_1,v_2)\in V: \, \theta(v)\ge 0 \,\rbrace,\\[6mm]
		&&\label{Am}(Au,v)_{V}= \int_a^{-l}{\cal F}_1 \Big(\frac{du_1}{dx}\Big)\frac{dv_1}{dx}\,dx+
		\int_l^{b}{\cal F}_2 \Big(\frac{du_2}{dx}\Big)\frac{dv_2}{dx}
		\, dx\\ [2mm]
		&&\qquad\qquad \forall\,u=(u_1,u_2),\,v=(v_1,v_2)\,\in V ,\nonumber\\[6mm]
		&&\label{jm}j(u,v)= -p(\theta(u))\theta(v)\quad \forall\,u=(u_1,u_2),\,v=(v_1,v_2)\,\in V ,\\[6mm]
		&&\label{fm}(f,v)_{V}=\int_a^{-l}f_{1}v_1\, dx+
		\int_l^{b}f_{2}v_2\, dx
		\quad \forall\,v=(v_1,v_2)\,\in V.
	\end{eqnarray}
	
	Note that  the definitions (\ref{Am}) and (\ref{fm}) are based on the Riesz  representation theorem. Moreover, here and below, we do not specify the dependence of various functions on the spatial variable $x$.

	With these preliminaries we are in a position to derive the variational formulation of Problem~$\cal P$. We assume in what follows that $u=(u_1,u_2)$ and $\sigma=(\sigma_1,\sigma_2)$ are sufficiently regular functions which satisfy (\ref{2.1})--(\ref{2.4}) and let $v=(v_1,v_2)\in K$. First, we perform an integration by parts  and use the equilibrium equation (\ref{2.2}) to see that
	\begin{eqnarray*}
		&&\int_a^{-l}\,\sigma_1\Big(\frac{dv_1}{dx}- \frac{du_1}{dx} \Big)\,dx=\int_a^{-l}\,
		f_1(v_1- u_1)\,dx\\ [2mm] 
		&&\qquad+\sigma_1(-l)(v_1(-l)-u_1(-l)) - \sigma_1(a)(v_1(a)-u_1(a)).
	\end{eqnarray*}

\noindent
	Next, since $v_1(a)=u_1(a) = 0$, we deduce that 
	\begin{equation}\label{xx}
		\int_a^{-l}\,\sigma_1\Big(\frac{dv_1}{dx}- \frac{du_1}{dx} \Big)\,dx=\int_a^{-l}	f_1(v_1- u_1)\,dx + \sigma_1(-l)(v_1(-l)-u_1(-l)).
	\end{equation}

\noindent
	A similar argument leads to equality
	\begin{equation}\label{xxn}
		\int_l^{b}\,\sigma_2\Big(\frac{dv_2}{dx}- \frac{du_2}{dx} \Big)\,dx=\int_l^b
		f_2(v_2-u_2)\,dx -\sigma_2(l)(v_2(l)-u_2(l)).
	\end{equation}
	Moreover, using notation (\ref{tm}) and equality $\sigma_1(-l)=\sigma_2(l)$ in (\ref{2.4}), we write 
	\begin{eqnarray*}
		&&\sigma_1(-l)(v_1(-l)-u_1(-l)) -\sigma_2(l)(v_2(l)-u_2(l))\\ [2mm]
		&&\quad = \sigma_1(-l)\big(v_1(-l)-v_2(l)-u_1(-l)+u_2(l)\big)
		\\ [2mm]
		&&\qquad = \sigma_1(-l)\big(\theta(u)-\theta(v)\big)
		\\ [2mm]&&\qquad 
		=\big(\sigma_1(-l)+p(\theta(u))\big)
		\big(\theta(u)-\theta(v)\big)-p(\theta(u))\big(\theta(u)-\theta(v)\big).
  \end{eqnarray*}
  Therefore,
  \begin{eqnarray}
   &&\label{zz}\sigma_1(-l)(v_1(-l)-u_1(-l)) -\sigma_2(l)(v_2(l)-u_2(l))\\ [2mm]   
&&\qquad\qquad=\big(\sigma_1(-l)+p(\theta(u))\big)\theta(u)-\big(\sigma_1(-l)+p(\theta(u))\big)\theta(v)\nonumber\\ [2mm]
		&&\qquad\qquad\qquad-p(\theta(u))\big(\theta(u)-\theta(v)\big).\nonumber
  \end{eqnarray}
Then, using (\ref{2.4})  and definition
	(\ref{Km}) of the set $K$ it is easy to see that
	\[\big(\sigma_1(-l)+p(\theta(u))\big)\theta(v)\le 0,\]
	\[\big(\sigma_1(-l)+p(\theta(u))\big)\theta(u)=0,\]
	hence \eqref{zz} implies that
	\begin{eqnarray}
		&&\label{yx}\sigma_1(-l)(v_1(-l)-u_1(-l)) -\sigma_2(l)(v_2(l)-u_2(l))\\ [2mm]
		&&\qquad\qquad\ge -p(\theta(u))\big(\theta(u)-\theta(v)\big).\nonumber
	\end{eqnarray}

	We now combine  (\ref{xx}), (\ref{xxn}), (\ref{yx}) and use the definitions (\ref{jm}), (\ref{fm}) to deduce that
	\begin{eqnarray}
		&&\label{Pb}\int_a^{-l}\,\sigma_1\Big(\frac{dv_1}{dx}- \frac{du_1}{dx} \Big)\,dx+\int_l^{b}\,\sigma_2\Big(\frac{dv_2}{dx}- \frac{du_2}{dx} \Big)\,dx
		\\ [3mm]
		&&\qquad
		+j(u, v) - j(u, u) \ge(f, v- u)_{V}\quad  \forall \ v\in K.\nonumber
	\end{eqnarray}
	We now substitute the constitutive laws  (\ref{2.1}), (\ref{2.1n}) in (\ref{Pb}) and use definition (\ref{Am}) of the operator $A$
	to obtain the following variational formulation of Problem~${\cal P}$ in term of displacements.

	\medskip
	\noindent{\bf Problem ${\cal P}_V$}. {\it Find  a displacement field $u$  such that
		the inequality below holds:}
	\begin{equation}\label{vari}
		u\in K,\qquad( Au, v-u)_{V} + j(u, v) - j(u, u) \ge(f, v- u)_{V}\quad \, \forall\,  v\in K.
	\end{equation}
	
	Our existence and uniqueness result in the study of Problem~${\cal P}_V$ is the following.
	
	\begin{Theorem}\label{t1}
		Assume $(\ref{F1})$--$(\ref{pc})$, $(\ref{sm})$ and $(\ref{fc})$. Then Problem  ${\cal P}_V$ has a unique solution $u\in K$.
	\end{Theorem}
	
	\noindent {\bf Proof}. We use Theorem \ref{t0a} on the space $X=V$.
	First, it is easy to see that the set $K$ given by (\ref{Km}) satisfies condition (\ref{Ka}). Next, we use  assumptions (\ref{F1}) and (\ref{F2}) to see that  the operator $A$ defined by (\ref{Am}) satisfies the inequalities
	\begin{eqnarray*}
		&&(Au-Av, u-v)_{V}\ge (m_1+m_2)\,\|u-v\|^2_{V}\qquad\forall\,u,
		v \in V, 
		\\[2mm]
		&&\|Au-Av\|_{V}\leq (M_1+M_2)\|u-v\|_{V} \qquad \forall\,u,\,v \in V.
	\end{eqnarray*}
This implies that condition \eqref{Aa} holds, too, with $m=m_1+m_2$.
	Next, we turn on the properties of the function $j$ defined by (\ref{jm}). Let $\eta\in V$. We note that the functional $j(\eta , \cdot)$ is an  affine continuous function and, therefore, is convex and lower semicontinuous. Consider now the elements $\eta=(\eta_{1},\eta_{2}),\,  \widetilde{\eta}=(\widetilde{\eta}_{1},\widetilde{\eta}_{2}),\,v=(v_{1},v_{2}), \, \widetilde{v}=(\widetilde{v}_{1},\widetilde{v}_{2})\in V$. Then, using definition (\ref{jm}), assumption  (\ref{pc})(b) and inequalities (\ref{injn}),  we find that 
	\begin{eqnarray*}
		&&\hspace{-3mm}j(\eta, \widetilde{v})-j(\eta, v)+j(\widetilde{\eta}, v)-j(\widetilde{\eta}, \widetilde{v})=\big(p(\theta(\eta))-p(\theta(\widetilde{\eta})\big)\big(\theta(v))-\theta(\widetilde{v}))\big)\\[2mm]
		&&\ \le L_p\big( |\eta_1(-l)-\widetilde{\eta}_1(-l)|+|\eta_2(l)-\widetilde{\eta}_2(l)|\big)\big( |v_1(-l)-\widetilde{v}_1(-l)|+|v_2(l)-\widetilde{v}_2(l)|\big)
		\nonumber\\[2mm]
		&&\ \ \le  L_pL\big(\|\eta_1-\widetilde{\eta}_1\|_{V_1}+\|\eta_2-\widetilde{\eta}_2\|_{V_2}\big)\big(\|v_1-\widetilde{v}_1\|_{V_1}+\|v_2-\widetilde{v}_2\|_{V_2}\big).\nonumber
	\end{eqnarray*}
	Next, using the inequalities,
	\begin{eqnarray*}
		&&\|\eta_1-\widetilde{\eta}_1\|_{V_1}+\|\eta_2-\widetilde{\eta}_2\|_{V_2}\le\sqrt{2}\, \|\eta-\widetilde{\eta}\|_{V},\\ [2mm]
		&&\|v_1-\widetilde{v}_1\|_{V_1}+\|v_2-\widetilde{v}_2\|_{V_2}\le\sqrt{2}\, \|v-\widetilde{v}\|_{V}
	\end{eqnarray*}
	we deduce that
    \begin{equation}
        \label{jab}j(\eta, \widetilde{v})-j(\eta, v)+j(\widetilde{\eta}, v)-j(\widetilde{\eta}, \widetilde{v}) 
        \le  2L_pL\,\|\eta-\widetilde{\eta}\|_{V}\|v-\widetilde{v}\|_{V}.
    \end{equation}
	
	Inequality (\ref{jab}) shows that $j$ satisfies condition (\ref{ja})	
	with $\alpha=2L_pL$.  Moreover, using  assumption (\ref{sm}) we deduce that the smallness condition  \eqref{small} holds, too. Theorem \ref{t1} is now a direct consequence of Theorem \ref{t0a}. \hfill$\Box$
	
	\medskip
	Once the displacement field $u=(u_1,u_2)$ is known, the stress field $\sigma=(\sigma_1,\sigma_2)$ can be easily obtained by using the constitutive laws (\ref{2.1}) and (\ref{2.1n}). A couple $(u, \sigma)$ with  $u=(u_1,u_2)$ and $\sigma=(\sigma_1,\sigma_2)$ which satisfies (\ref{2.1}),  (\ref{2.1n}) and (\ref{vari}) is called a weak solution to contact Problem~${\cal P}$. We conclude by Theorem~\ref{t1} that Problem~${\cal P}$ has a unique weak solution.

	\section{Convergence results}\label{s4}
	\setcounter{equation}0

	In this section we apply Theorem \ref{t1a} in the study of Problem ${\cal P}_V$. To this end, everywhere below we
	assume that $(\ref{F1})$--$(\ref{pc})$, $(\ref{sm})$ and $(\ref{fc})$  hold, even if we do not mention it explicitely.
	Consider now  a function $q$ which satisfies the following condition.
	\begin{equation}
		\hspace{10mm}	\left\{\begin{array}{ll} {\rm (a)\ }
			q\colon \mathbb{R}\to\mathbb{R}.\\[2mm]
			{\rm (b)\ There\ exists}\, L_q>0 {\rm\ such\ that\ }\\
			{}\ \qquad  |q(r_1)-q(r_2)|
			\le L_q\,|r_1-r_2|\quad\forall\,
			r_1,\,r_2\in \mathbb{R}.\qquad\qquad\\[2mm]
			{\rm (c)\ }(q(r_1)-q(r_2))(r_1-r_2)\le 0\quad\forall\,
			r_1,\,r_2\in \mathbb{R}.\\[2mm]
			{\rm (c)\ }q(r)\ge 0 \ \ {\rm if\ }\ r \le2l \ \ {\rm and}\ \ q(r)\le 0\ \ {\rm if\ }\ r \ge 2l.
		\end{array}\right.
		\label{qc}
	\end{equation}
	Then, we use the Riesz representation theorem to define the operator  $G\colon V \longrightarrow V$ by equality
	\begin{equation}\label{Gp}
		(Gu, v)_{V} = q(\theta(u))v_1(-l)-q(\theta(u))v_2(l)\quad \forall\,u=(u_1,u_2),\,v=(v_1,v_2)\,\in V.
	\end{equation}
We note that
\begin{equation}\label{Gpn}
	(Gu, v)_{V} = q(\theta(u))(2l-\theta(v))\quad \forall\,u=(u_1,u_2),\,v=(v_1,v_2)\,\in V
\end{equation}
and, combining this equality with the properties of the functions $q$ and $\theta$, we deduce that
\begin{equation}\label{GG}
	\mbox{$G\colon V\to V$  is a  monotone Lipschitz continuous operator.}
\end{equation}
	
Consider now a sequence $\{\lambda_n\}\subset\real$ such that $\lambda_n>0$. For each $n\in\mathbb{N}$
denote by $\cP_V^n$  the following penalty version of Problem $\cP_V$.

\medskip\noindent{\bf Problem}  ${\cal P}_V^n$. {\it Find $u_n$ such that}
\begin{eqnarray}\label{1m}
	&&u_n\in K,\quad (A u_n, v - u_n)_V +
	\frac{1}{\lambda_n} (Gu_n, v - u_n)_V +
	j(u_n,v) - j(u_n,u_n)\label{2m}\\[2mm]
	&&\qquad \qquad\qquad
	\ge (f, v - u_n)_V
	\qquad\forall\, v \in K.\nonumber
	\nonumber
\end{eqnarray}

Then, arguments similar to those used in the proof of Theorem \ref{t1}, based on properties \eqref{GG} and Theorem \ref{t0a}, imply that
Problem $\cP_V^n$ has a unique solution $u_n\in K$, for each $n\in\mathbb{N}$.

We now define the sets $K'$, $K''$ and $K'''$  by equalities: 
\begin{eqnarray}
	&&\label{Kp}K' = \lbrace v=(v_1,v_2)\in V: \, \theta(v)\ge 2l \,\rbrace,\\[2mm]
	&&\label{Kpp}K'' = \lbrace v=(v_1,v_2)\in V: \, 0\le \theta(v)\le 2l \,\rbrace,\\[2mm]
	&&\label{Kppp}K''' = \lbrace v=(v_1,v_2)\in V: \, \theta(v)=2l \,\rbrace.
\end{eqnarray}
We associate to these sets the following inequality problems.

\medskip
\noindent{\bf Problem ${\cal P}'_V$}. {\it Find a displacement field $u'$  such that
	the inequality below holds:}
\begin{equation}\label{vp}
	u'\in K',\qquad( Au', v-u')_{V} + j(u', v) - j(u', u') \ge(f, v- u')_{V}\quad \, \forall \, v\in K'.
\end{equation}

\medskip
\noindent{\bf Problem ${\cal P}''_V$}. {\it Find a displacement field $u''$  such that
	the inequality below holds:}
\begin{equation}\label{vpp}
	u''\in K'',\qquad( Au'', v-u'')_{V} + j(u'', v) - j(u'', u'') \ge(f, v- u'')_{V}\quad \, \forall\,  v\in K''.
\end{equation}

\medskip
\noindent{\bf Problem ${\cal P}'''_V$}. {\it Find a displacement field $u'''$  such that
	the inequality below holds:}
\begin{equation}\label{vppp}
	u'''\in K''',\quad( Au''', v-u''')_{V} + j(u''', v) - j(u''', u''') \ge(f, v- u''')_{V}\quad \forall \, v\in K'''.
\end{equation}

\medskip
Note that these problems are similar to Problem $\cP_V$, the only difference arising in the fact that here the set of constraints $K$ was successively replaced by the sets $K'$, $K''$ and $K'''$, respectively. Nevertheless, since these sets are nonempty convex closed subsets of $V$, we deduce that each of Problems 
${\cal P}'_V$, ${\cal P}''_V$ and ${\cal P}'''_V$
has a unique solution, denoted in what follows by $u'$, $u''$ and $u'''$, respectively.	

\medskip

Next, we assume that the function $q$ successively satisfies one of the following additional conditions.
\begin{equation}\label{qqq}
\left\{\begin{array}{ll}
{\rm (i)}\quad \,\,q(r)=0\ \ \mbox{if and only if}\ \ r\ge 2l.\\ [2mm]
{\rm (ii)}\quad\, q(r)=0\ \ \mbox{if and only if}\ \ r\le 2l.\\ [2mm]
{\rm (iii)}\quad q(r)=0\ \ \mbox{if and only if}\ \ r=2l
\end{array}\right.
\end{equation}
Note that examples of functions $q$ which satisfy the above conditions can be easily constructed. For instance, 
an example of a function $q$ which satisfies (\ref{qc}) and (\ref{qqq})(i) is given by
\begin{equation}
	q(r)=
	\left\{\begin{array}{ll}
		2l-r\qquad {\rm if}\quad r < 2l,\\[2mm]
		0\qquad \qquad{\rm if}\quad r\ge 2l.\end{array}\right.
	\label{qq}
\end{equation}
Moreover, an example of a function $q$ which satisfies (\ref{qc}) and (\ref{qqq})(ii) is given by
\begin{equation}
	q(r)=
	\left\{\begin{array}{ll}
		0\qquad\quad {\rm if}\quad r < 2l,\\[2mm]
		2l-r\quad{\rm if}\quad r\ge 2l.\end{array}\right.
	\label{qqn}
\end{equation}	

\medskip

We now state and prove our main result in this section.

\begin{Theorem}\label{t2}
	Under assumptions $(\ref{F1})$--$(\ref{pc})$, $(\ref{sm})$, $(\ref{fc})$, $(\ref{qc})$, $(\ref{za5})$ and $(\ref{za6})$, the following statements hold.

\smallskip
{\rm a)} If the function $q$ satisfies condition \eqref{qqq} {\rm (i)}, then $u_n\to u'$ in $V$.	

\smallskip
{\rm b)} If the function $q$ satisfies condition  \eqref{qqq} {\rm (ii)}, then $u_n\to u''$ in $V$.

\smallskip
{\rm c)} If the function $q$ satisfies condition \eqref{qqq} {\rm (iii)}, then $u_n\to u'''$ in $V$.

\end{Theorem}

	\noindent {\bf Proof}. 
	a) We apply Theorem \ref{t1a} with $ X= V$ and $K^*=K'$. Since the rest of the conditions in this theorem are obviously satisfied, we only have to check that the operator $G$ satisfies conditions (\ref{za4}), that is
	\begin{equation}
	\label{za4n}\left\{\begin{array}{ll}
		\mbox{\rm (a)}
		\quad(Gu,v-u)_V\le 0\qquad\forall\, u\in K,\ v\in K',\\ [3mm]
		\mbox{\rm (b)} \quad u\in K,\quad (Gu,v-u)_V=0\quad\forall\,v\in K'\ \ \Longrightarrow\ \  u\in K'.
	\end{array}\right.\\ [2mm]
	\end{equation}

Let $u\in K$ and $v\in K'$. Then, using \eqref{Kp} and \eqref{qqq}(i) we find that
$q(\theta(v))=0$.  Therefore, using
\eqref{Gpn}, we have
\[
(Gu, v-u)_{V} = q(\theta(u))(\theta(u)-\theta(v))=
 (q(\theta(u))-q(\theta(v)))(\theta(u)-\theta(v))\]
 and, invoking assumption (\ref{qc})(c), we deduce that $(Gu,v-u)_V\le 0$. 
 Assume now that $(Gu,v-u)_V=0$ for all $v\in K'$. Then 
 \[q(\theta(u))(\theta(u)-\theta(v))=0\qquad\forall\,v\in K'\] 
 which implies that either $\theta(u)=\theta(v)$ or $q(\theta(u))=0$.
 We use definition (\ref{Kp}) and assumption \eqref{qqq}(i) to deduce that in both cases 
$\theta(u)\ge 2l$ and, therefore, $u\in K'$.
It follows from above that condition \eqref{za4n} is satisfied, which concludes the proof of this part.
	
	\medskip
	
	b), c) We apply Theorem \ref{t1a} with $ X= V$ and  $K^*=K''$, $K^*=K'''$, respectively.  The arguments are similar to those use in the part a) of the proof and, therefore, we skip the details. \hfill$\Box$
	
\medskip We end this section with the following mechanical interpretation of Theorem \ref{t2}.

\medskip
a)	 First,  Problem $\cP_V^n$ represents the variational formulation of  a version of Problem~$\cP$ in which  condition (\ref{2.4}) is replaced by  condition
\begin{eqnarray}
	&&\left\{ 
	\begin{array}{ll}
		\theta(u)\ge 0,\\[2mm]
		\sigma_1(-l)=\sigma_2(l)=-p(\theta(u))-\frac{1}{\lambda_n}\, q(\theta(u))\quad{\rm if}\quad \theta(u)>0,\\[2mm]
		\sigma_1(-l)=\sigma_2(l)
		\le-p(\theta(u))-\frac{1}{\lambda_n}\, q(\theta(u))\quad{\rm if}\quad \theta(u)=0.
		\label{2.4a}\end{array}\right.
\end{eqnarray} 
It models the equilibrium of the spring-rods system in the case when the behavior of the spring is described with the function $p+\frac{1}{\lambda_n}\, q$, in which $\frac{1}{\lambda_n}$ can be interpreted as an additional stiffness coefficient.

\medskip
b)	 Second,  Problem $\cP_V'$ represents the variational formulation of a version of Problem $\cP$ in which the condition (\ref{2.4}) is replaced by the condition
\begin{eqnarray}
	&&\left\{ 
	\begin{array}{ll}
		\theta(u)\ge 2l,\\[2mm]
		\sigma_1(-l)=\sigma_2(l)=-p(\theta(u))\quad{\rm if}\quad \theta(u)>2l,\\[2mm]
		\sigma_1(-l)=\sigma_2(l)
		\le-p(\theta(u))\quad{\rm if}\quad \theta(u)=2l.
		\label{2.4b}\end{array}\right.
\end{eqnarray} 
It models the equilibrium of the spring-rods system in the limit case when the spring behaves like a rigid in compression but it behaves elastically in extension. 
Theorem~\ref{t2}~a) shows that, whichever are the applied forces, the weak solution of this problem can be approached by the weak solution of Problem $\cP$ with a large stiffness in compression, provided that
the stiffness in extension does not change.

\medskip
c)	 Third,  Problem $\cP_V''$ represents the variational formulation of a version of Problem $\cP$ in which the condition (\ref{2.4}) is replaced by the condition
\begin{eqnarray}
	&&\left\{ 
	\begin{array}{ll}
		0\le \theta(u)\le 2l,\\[2mm]
		\sigma_1(-l)=\sigma_2(l)=-p(\theta(u))\, \quad{\rm if}\quad 0\le \theta(u)<2l,\\[2mm]
		\sigma_1(-l)=\sigma_2(l)
		\le-p(\theta(u))\quad{\rm if}\quad \theta(u)=2l.
		\label{2.4c}\end{array}\right.
\end{eqnarray} 

It models the equilibrium of the spring-rods system in the limit case when the spring behaves as rigid in extension, but it behaves elastically in compression. 
Theorem~\ref{t2}~b) shows that, whichever are the applied forces, the weak solution of this problem can be approached by the weak solution of Problem~$\cP$ with a large stiffness in extension, provided that the behavior of the spring in compression does not change.

\medskip
d)	 Finally,  Problem $\cP_V'''$ represents the variational formulation of a version of Problem~$\cP$ in which the condition (\ref{2.4}) is replaced by the condition
\begin{eqnarray}
	&&\left\{ 
	\begin{array}{ll}
		\theta(u)=2l,\\[2mm]
		\sigma_1(-l)=\sigma_2(l)=0.
		\label{2.4d}\end{array}\right.
\end{eqnarray} 
It models the equilibrium of the spring-rods system in the limit case when the spring behaves as rigid in both extension and compression. 
Theorem~\ref{t2}~c) shows that, whichever are the applied forces, the weak solution of this problem can be approached by the weak solution of Problem~$\cP$ with a large stiffness in both extension and compression.

	\section{Numerical simulations}\label{s5}
	\setcounter{equation}0

In this section we present numerical simulations in the study of Problem~${\cal P}_V$},  which illustrate the theoretical results obtained in previous sections. For simplicity, we restrict ourselves to the case when the elasticity operators ${\cal F}_i$ (with $i=1,2$)  are linear, do not depend on the spatial variable $x$  and the function~$p$ is given by \eqref{pp}.
Therefore, we assume that ${\cal F}_i \colon \mathbb{R}\to \mathbb{R}$ are given by
\begin{equation*}
{\cal F}_1(x,r)= E_1 r\quad\forall x\in(a,-l),\ r\in  \mathbb{R},\quad {\cal F}_2(x,r)= E_2 r\quad\forall x\in(l,b),\ r\in  \mathbb{R},
\end{equation*}
where $E_1$ and $E_2$ represent the Young moduli of the rods.
Let $\hat p \colon \mathbb{R}\to \mathbb{R}$ be the function defined by
\begin{equation}\label{cone}
\displaystyle \hat p(r)=-\int_{2l}^rp(s)\,ds=
    \left\{ 
	\begin{array}{ll}
            \frac{k_1}{2}(r-2l)^2 &\ \ \text{if }\ \  r< 2l,\\ [2mm]
	\frac{k_2}{2}(r-2l)^2 &\ \  \text{if }\ \  r\geq 2l.
        \end{array} \right.
\end{equation}
Since $p$ is a decreasing function it follows that $\hat p$ is a convex function and, therefore, using the subgradient inequality, we deduce that
\begin{equation}\label{pp1}
\hat p(s)-\hat p(r)\ge -p(r)(s-r)\qquad\forall\, r,\, s\in\mathbb{R}.
\end{equation}
Therefore, using \eqref{jm}, \eqref{pp1}  and notation 
\begin{equation}
\varphi(v)=\hat p(\theta(v))\qquad\forall\, v\in V
\end{equation}
we deduce that 
\[\varphi(v)-\varphi(u)\ge j(u,v)-j(u,u)\qquad\forall\, u,\, v \in V
\]which shows that, if $u$ satisfies the inequality 
 (\ref{vari}), then
 \begin{equation}\label{vi}
u\in K,\qquad (Au, v-u)_{V} + \varphi(v) - \varphi(u) \ge(f, v- u)_{V}\quad \, \forall\,  v\in K.
 \end{equation}
 Moreover, since a well-known existence and uniqueness result implies that inequality (\ref{vi}) has a unique solution, we deduce that 
 inequalities~(\ref{vari})  and (\ref{vi}) are equivalent.
 
Note also that assumptions (\ref{cone}) imply that
the operator $A$  is linear and symmetric. Therefore,  using  notation
\begin{equation}
\label{minimization-problem}
F(v) = \frac{1}{2}\,(Av, v)_{V} + \varphi(v) - (f, v)_{V}\qquad\forall\, v\in V
\end{equation}
combined with a standard argument, we find that 
Problem $\cP_V$ is equivalent with the following minimization problem.

\medskip
\noindent{\bf Problem $\hat{\cal P}_V$}. {\it Find  a displacement field $u$  such that
	the inequality below holds:}
\begin{equation}\label{min}
u\in K,\qquad F(u)\le F(v)\qquad\forall\, v\in K.
\end{equation}

\medskip
It follows from above that  the numerical approximation of Problem $\cP_V$ could be carried out by using  the numerical approximation of its equivalent formulation $\hat{\cP}_V$.  We use this idea and, to perform the numerical approximation, we use arguments based on the Finite Element Method (FEM). Since the topic is standard, we skip the details. Nevertheless, we mention the following:
we used continuous piecewise affine finite elements in order to approach the solution of problem (\ref{min});
the code created to perform these simulations can be found in the repository https://github.com/KOS-UJ/Spring-Rods-System-Approximation; 
the approximate solution is the vector of displacement of discrete points in the rods, which minimizes the functional~$F$ given by~(\ref{minimization-problem}), approximated with FEM; 
the minimization problem \eqref{min} is solved using the Sequential Least Squares Programming implementation from the SciPy package for scientific computation \cite{scipy}.

In the numerical experiments below, we use the following data:
\[a=-1,\quad b=1,\quad l=0.5,\quad  E_1=E_2=1.\]
This choice makes the system symmetric in the reference configuration, when the forces are not applied. 
Note that, here and below, for simplicity, we do not indicate the units associated to the data and unknowns. In addition, as already mentioned, we consider a nonlinear spring whose behavior is described by the function $p$ given by \eqref{pp}. 
 Since in the equilibrium state a spring cannot be simultaneously compressed and extended, we 
 remark that a single simulation depends only on one stiffness parameter $k_1$ or $k_2$, and never on both. Finally, we assume that the densities of body forces are constant and will be described below. Therefore, it is easy to check that the above data satisfy assumptions \eqref{F1}--\eqref{pc} and \eqref{fc}, and to ensure the smallness assumption~\eqref{sm}, the spring stiffness coefficients are chosen such that $k_i\in (0,2)$, $i=1,2$.

\begin{figure}
        \centering
        \begin{subfigure}[b]{0.49\textwidth}
    	    \centering
    	    \includegraphics[width=\textwidth]{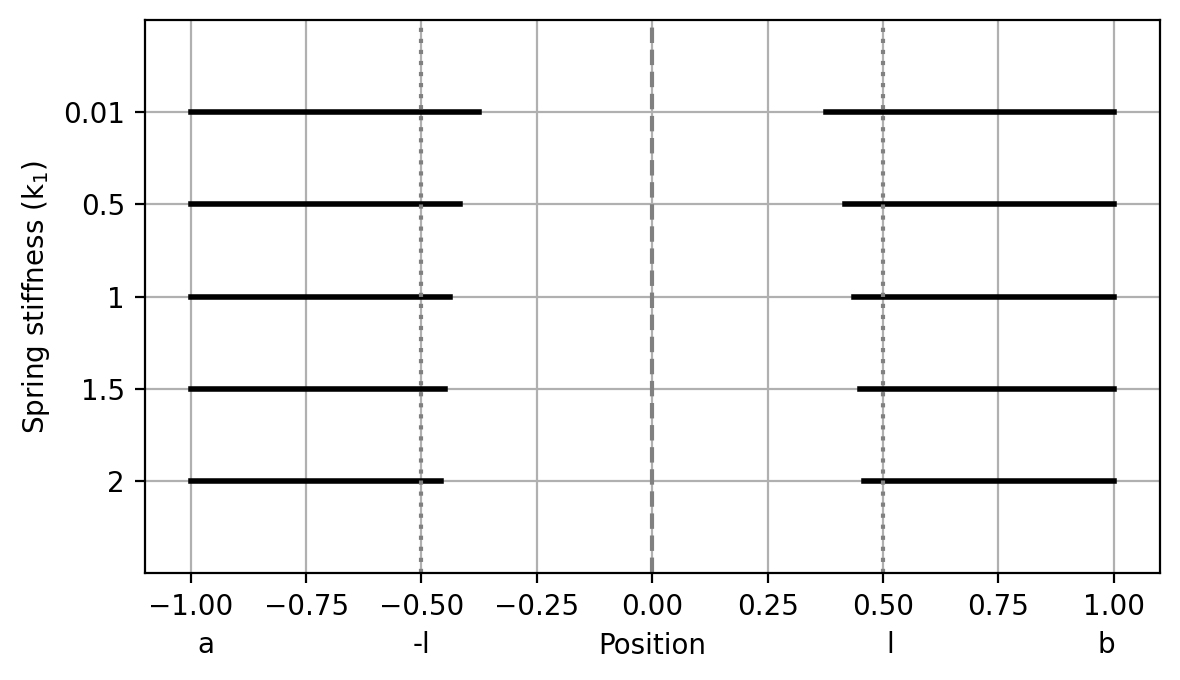}
    	    \caption{Equilibrium configuration.}
    	    \label{fig:experiment1_stress}
        \end{subfigure}
        \begin{subfigure}[b]{0.49\textwidth}
    	    \centering
    	    \includegraphics[width=\textwidth]{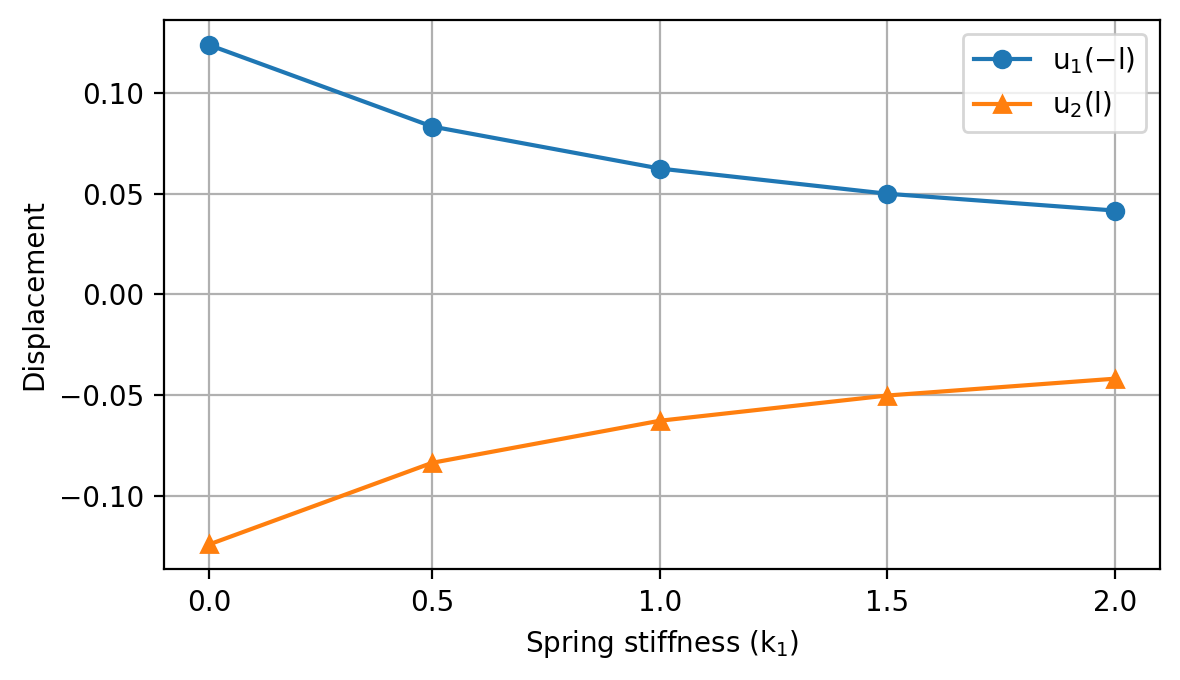}
    	    \caption{Displacement of the ends of rods.}
    	    \label{fig:experiment1_ends_displacements}
        \end{subfigure}
        \begin{subfigure}[b]{0.49\textwidth}
    	    \centering
    	    \includegraphics[width=\textwidth]{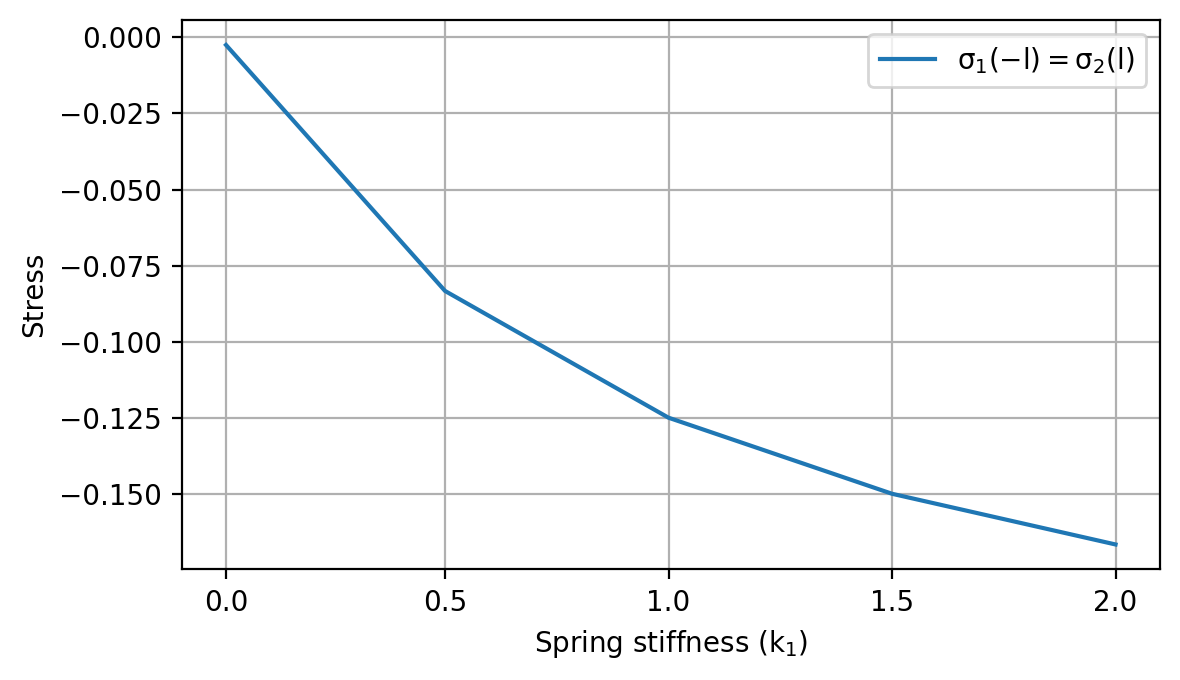}
    	    \caption{Stress on the ends of rods.}
    	    \label{fig:experiment1_ends_stress}
        \end{subfigure}
        
        \caption{Numerical results for Experiment \ref{e1}.}
        \label{fig:experiment1_figures}
    \end{figure}

\begin{figure}\label{fig2}
        \centering
        \begin{subfigure}[b]{0.49\textwidth}
    	    \centering
    	    \includegraphics[width=\textwidth]{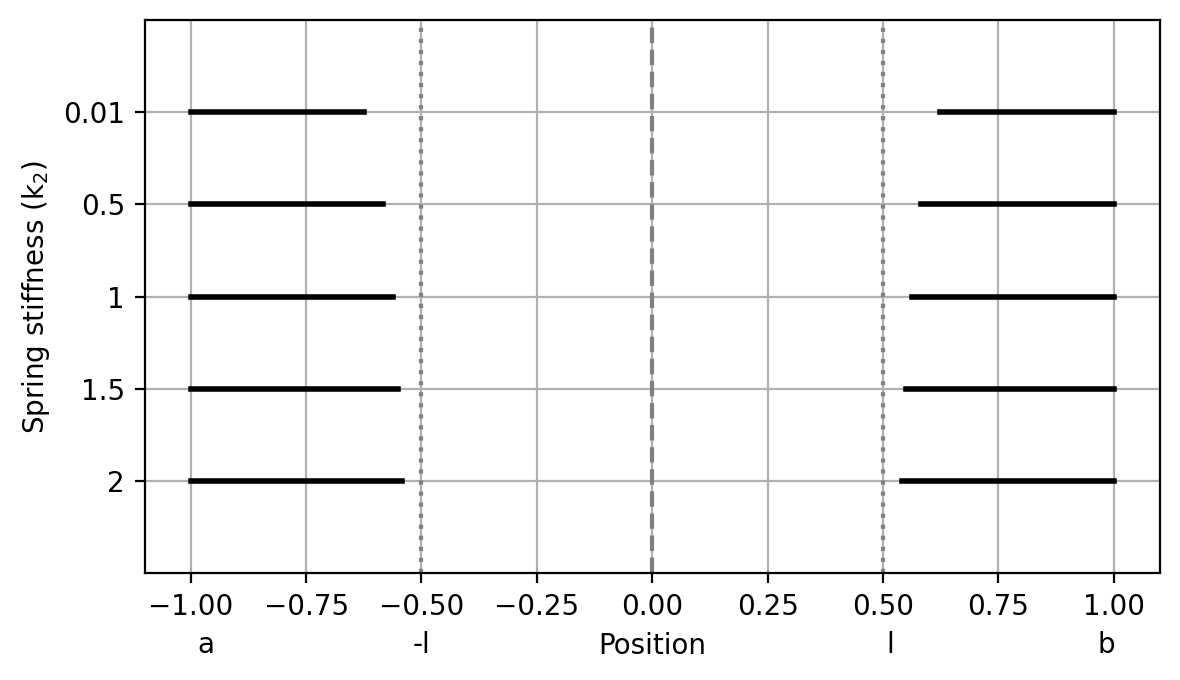}
    	    \caption{Equilibrium configuration.}
    	    \label{fig:experiment2_stress}
        \end{subfigure}
        \begin{subfigure}[b]{0.49\textwidth}
    	    \centering
    	    \includegraphics[width=\textwidth]{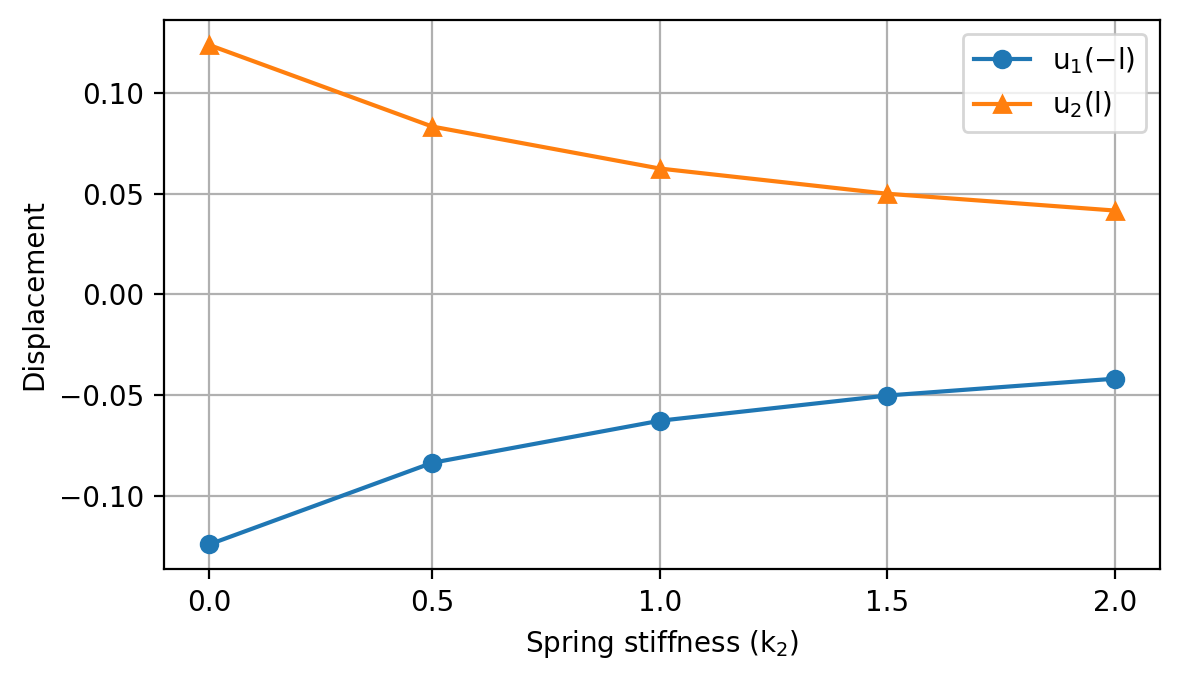}
    	    \caption{Displacement of the ends of rods.}
    	   \label{fig:experiment2_ends_displacements}
        \end{subfigure}
        \begin{subfigure}[b]{0.49\textwidth}
    	    \centering
    	    \includegraphics[width=\textwidth]{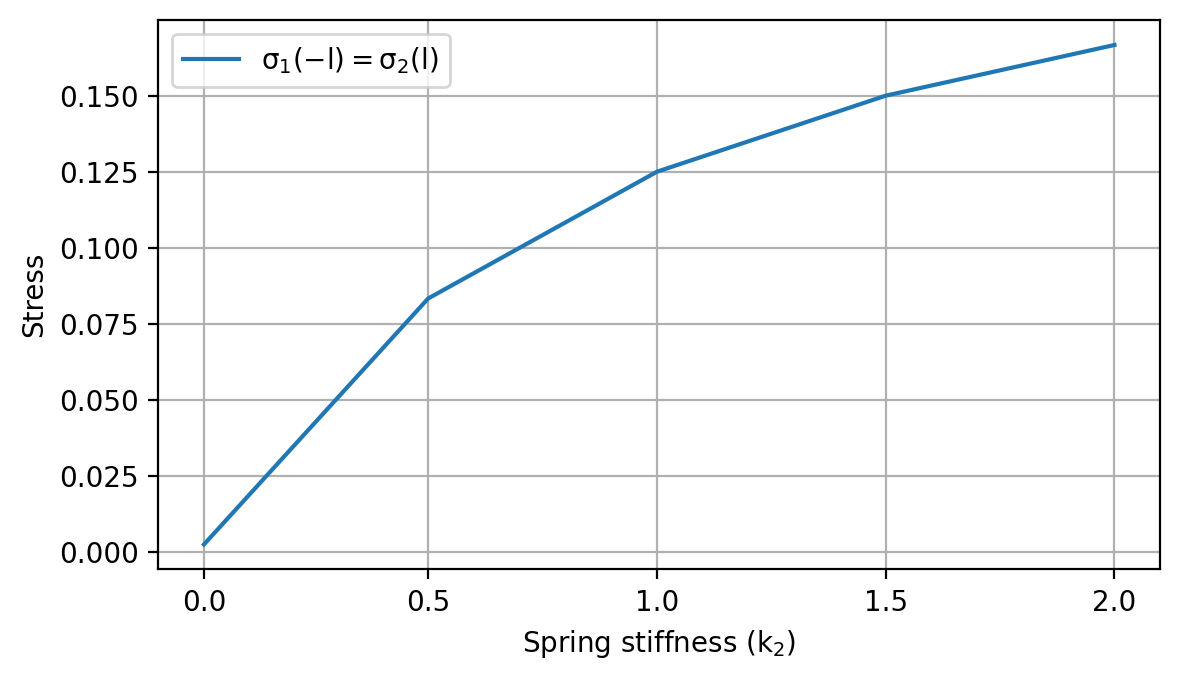}
    	    \caption{Stress on the ends of rods.}
    	    \label{fig:experiment2_ends_stress}
        \end{subfigure}
        
        \caption{Numerical results for Experiment \ref{e2}.} \label{fig:experiment2_figures}
    \end{figure}

\begin{figure}\label{fig3}
  \centering
        \begin{subfigure}[b]{0.49\textwidth}
    	    \centering
    	    \includegraphics[width=\textwidth]{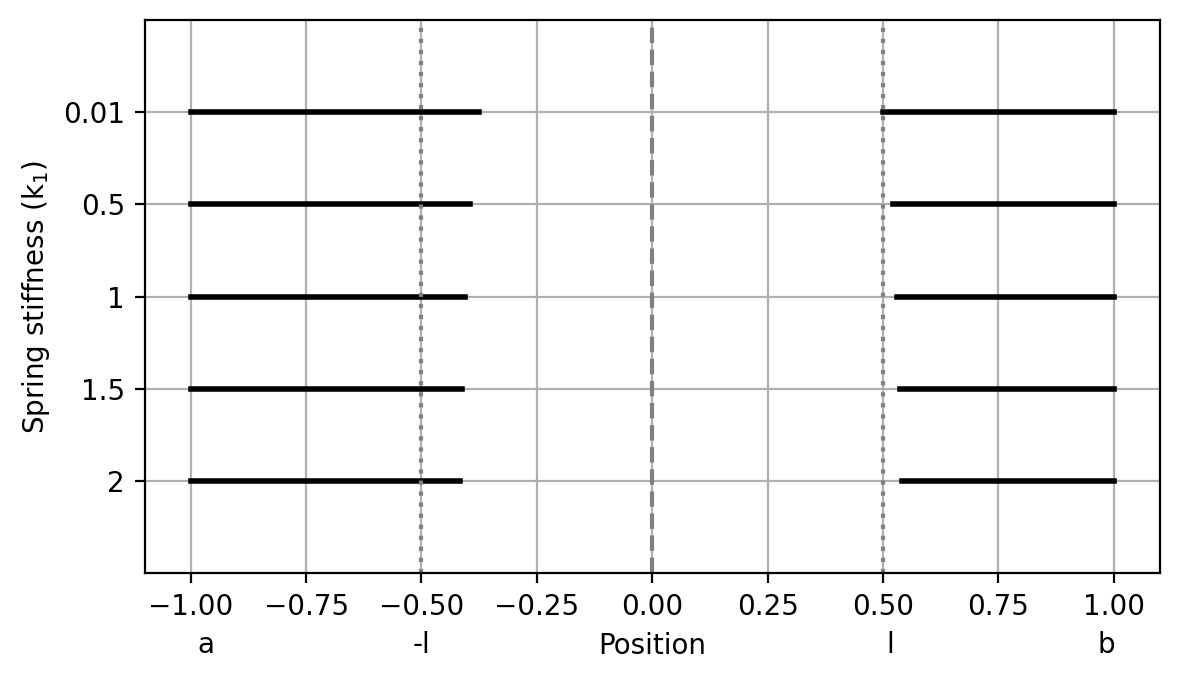}
    	    \caption{Equilibrium configuration.}
    	    \label{fig:experiment3_stress}
        \end{subfigure}
        \begin{subfigure}[b]{0.49\textwidth}
    	    \centering
    	    \includegraphics[width=\textwidth]{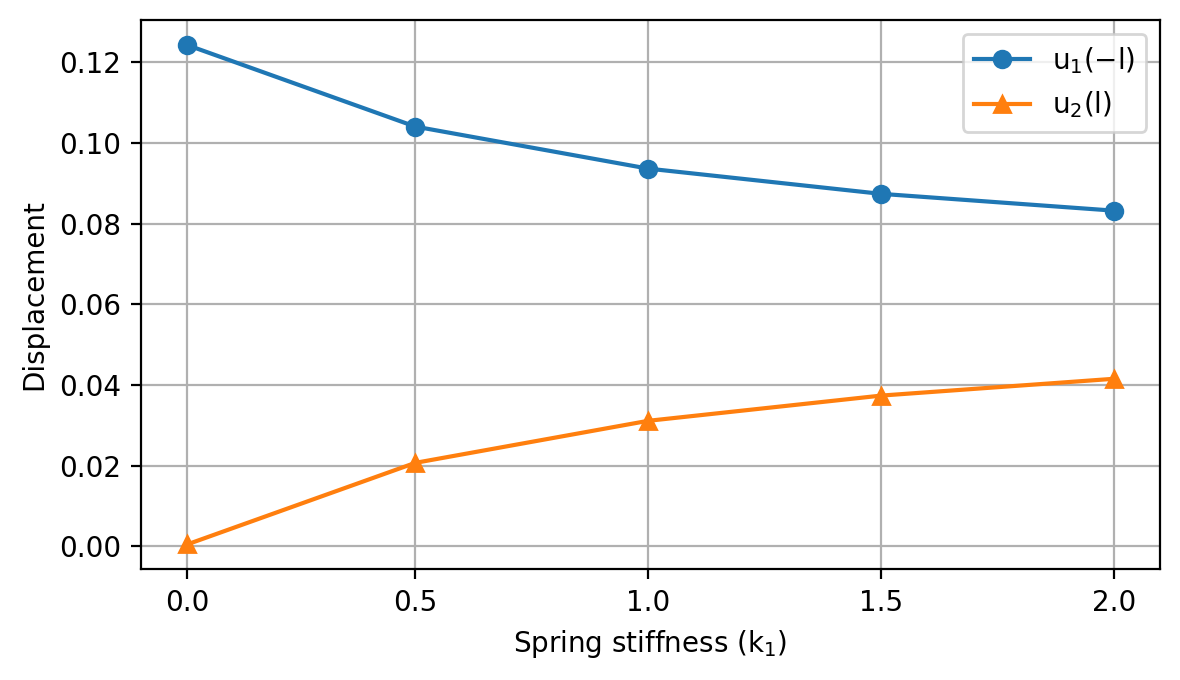}
    	    \caption{Displacement of the ends of rods.}
    	    \label{fig:experiment3_ends_displacements}
        \end{subfigure}
        \begin{subfigure}[b]{0.49\textwidth}
    	    \centering
    	    \includegraphics[width=\textwidth]{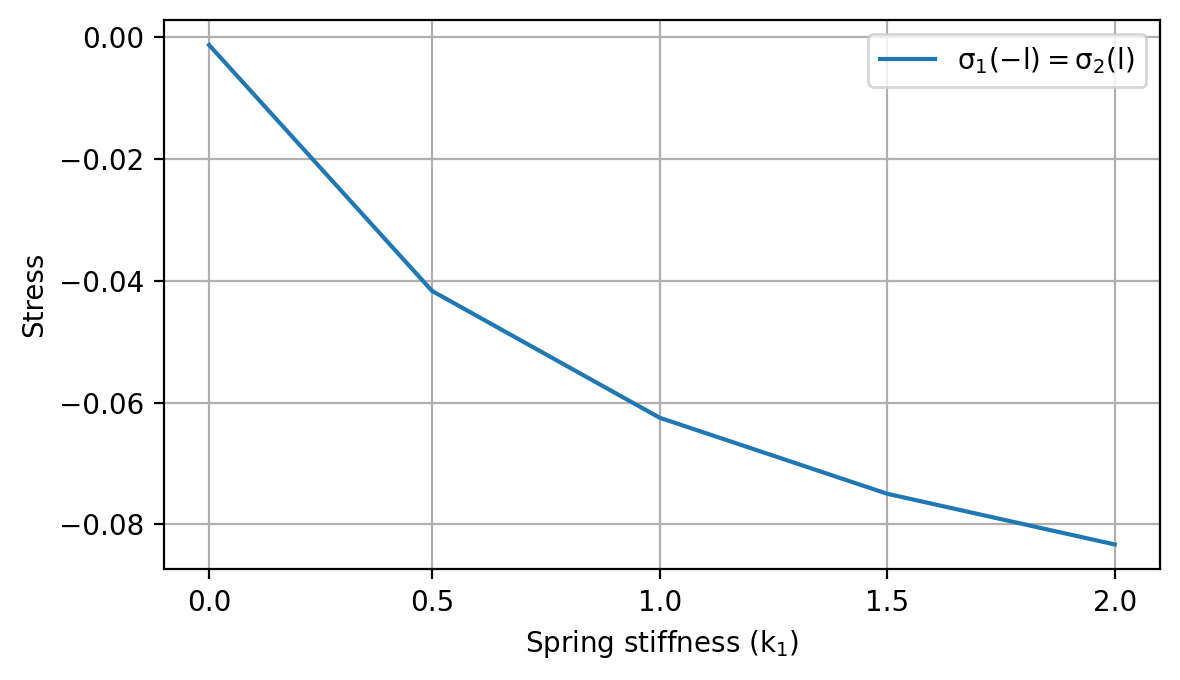}
    	    \caption{Stress on the ends of rods.}
    	    \label{fig:experiment3_ends_stress}
        \end{subfigure}
        
        \caption{Numerical results for Experiment \ref{e3}. }
        \label{fig:experiment3_figures}
    \end{figure}

\begin{figure}\label{fig5}
        \centering
        \begin{subfigure}[b]{0.49\textwidth}
    	    \centering
    	    \includegraphics[width=\textwidth]{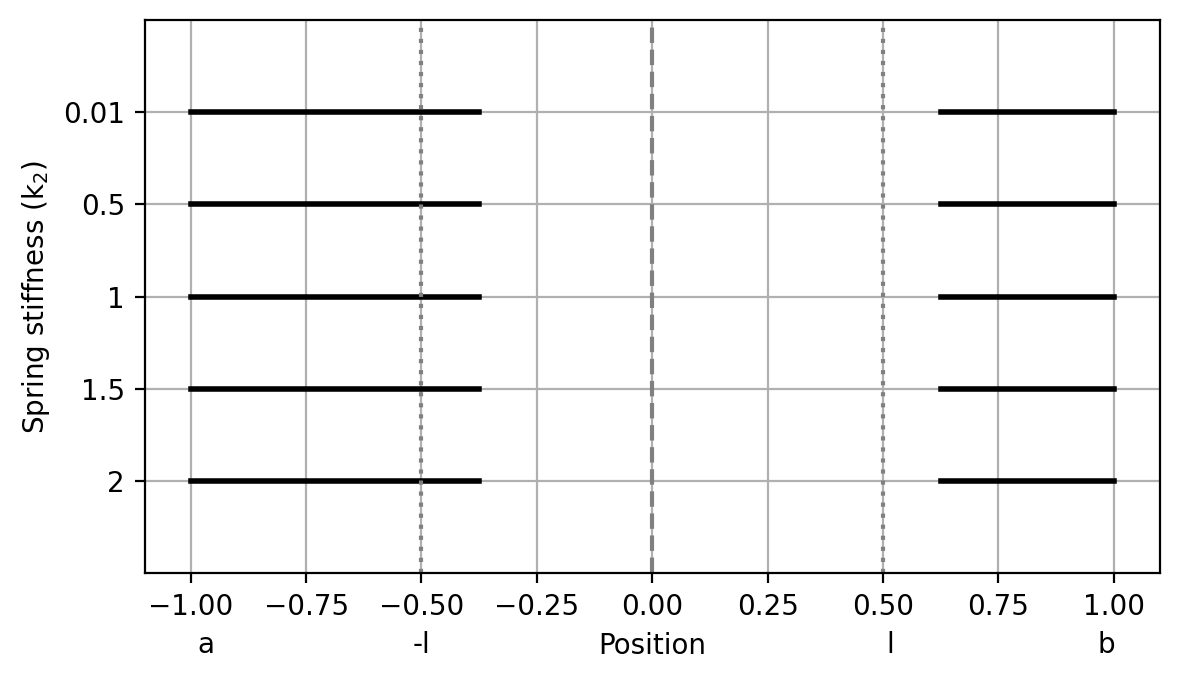}
    	    \caption{Equilibrium configuration.}
    	    \label{fig:experiment5_stress}
        \end{subfigure}
        \begin{subfigure}[b]{0.49\textwidth}
    	    \centering
    	    \includegraphics[width=\textwidth]{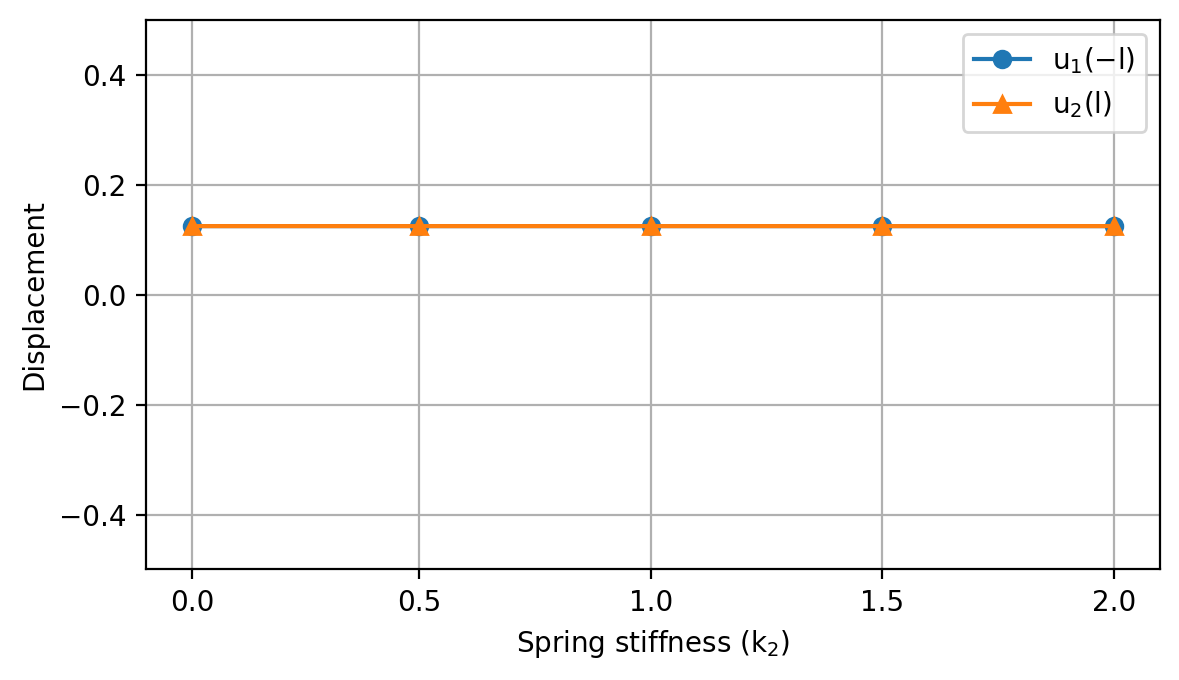}
    	    \caption{Displacement of the ends of rods.}
    	    \label{fig:experiment5_ends_displacements}
        \end{subfigure}
        \begin{subfigure}[b]{0.49\textwidth}
    	    \centering
    	    \includegraphics[width=\textwidth]{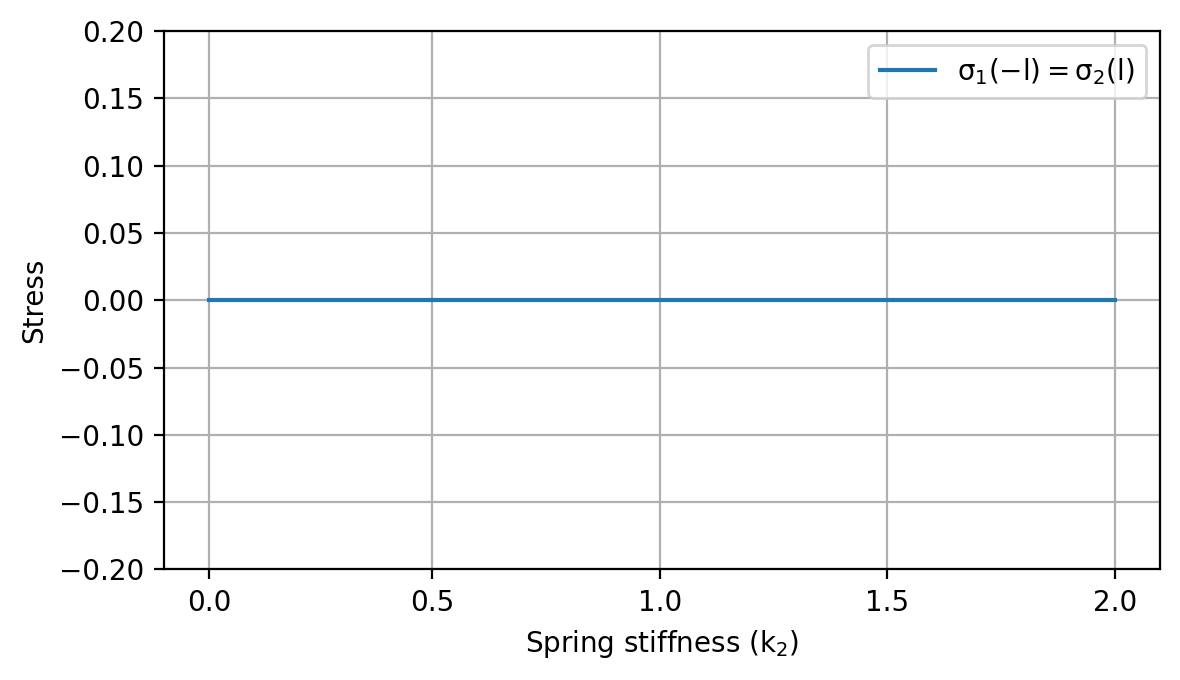}
    	    \caption{Stress on the ends of rods.}
    	    \label{fig:experiment5_ends_stress}
        \end{subfigure}
        
        \caption{Numerical results for Experiment \ref{e5}.}
        \label{fig:experiment5_figures}
    \end{figure}

\begin{figure}\label{fig5.5}
	\centering
	\begin{subfigure}[b]{0.49\textwidth}
		\centering
		\includegraphics[width=\textwidth]{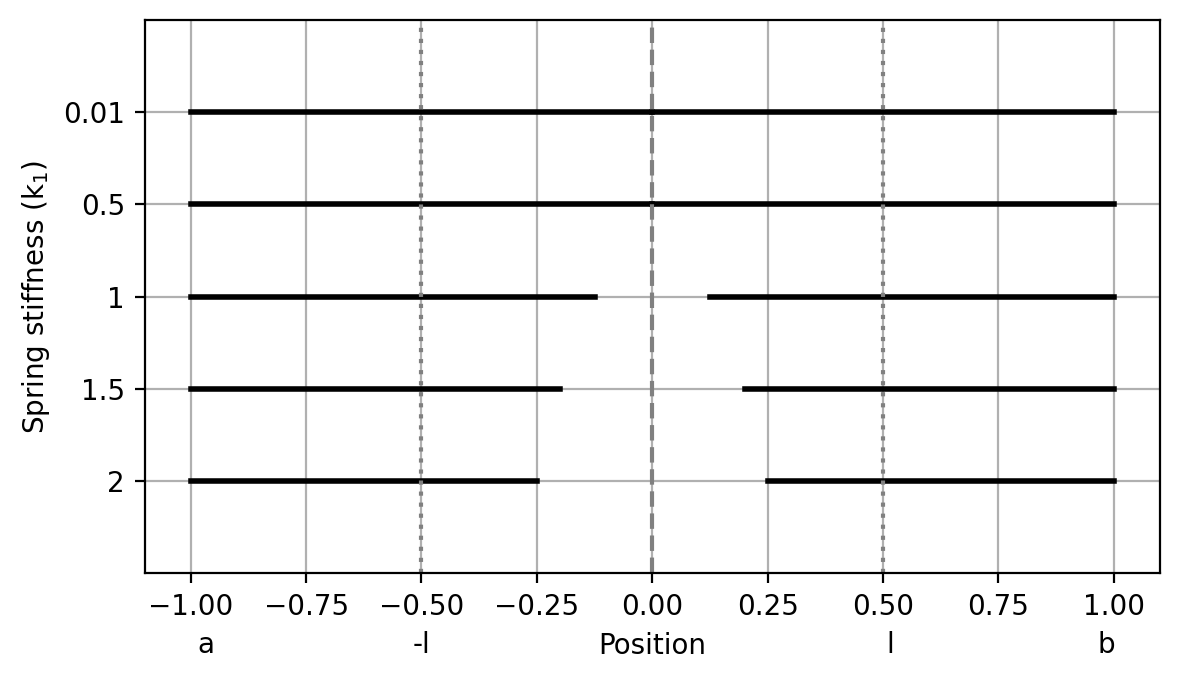}
		\caption{Equilibrium configuration.}
		\label{fig:experiment5.5_stress}
	\end{subfigure}
	\begin{subfigure}[b]{0.49\textwidth}
		\centering
		\includegraphics[width=\textwidth]{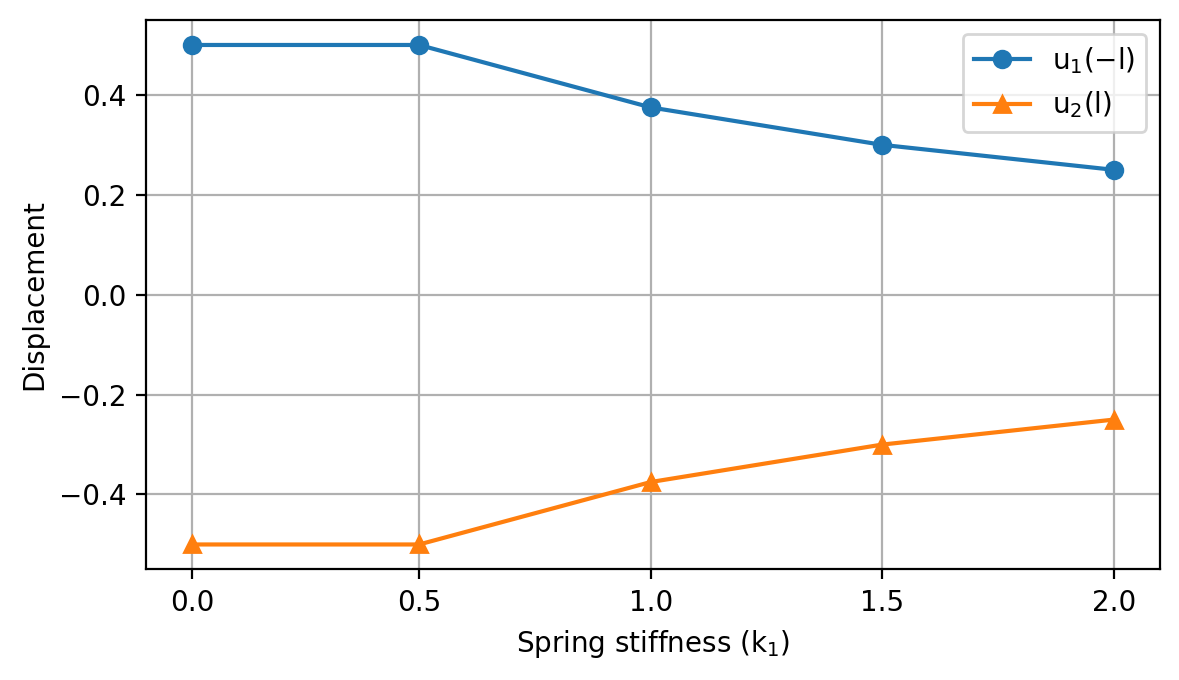}
		\caption{Displacement of the ends of rods.}
		\label{fig:experiment5.5_ends_displacements}
	\end{subfigure}
	\begin{subfigure}[b]{0.49\textwidth}
		\centering
		\includegraphics[width=\textwidth]{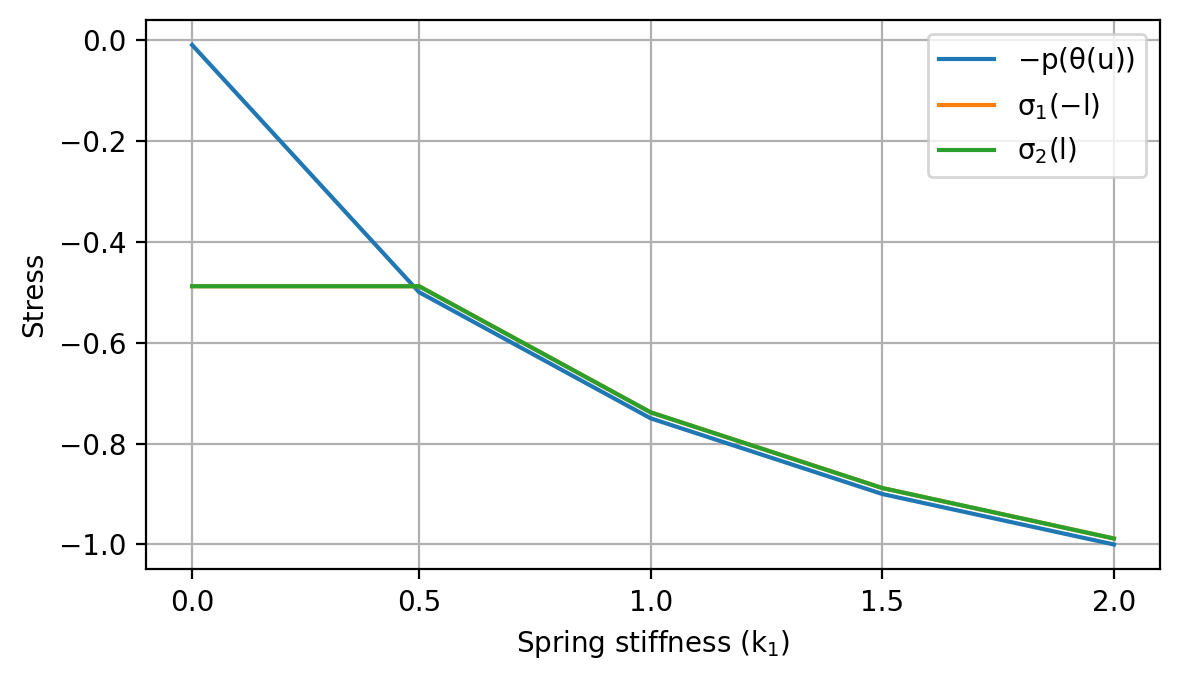}
		\caption{Stress on the ends of rods.}
		\label{fig:experiment5.5_ends_stress}
	\end{subfigure}
	
	\caption{Numerical results for Experiment \ref{e5.5}.}
	\label{fig:experiment5.5_figures}
	\end{figure}

\begin{figure}\label{fig6}
    	\vspace{-13mm}
        \centering
        \begin{subfigure}[b]{0.49\textwidth}
    	    \centering
    	    \includegraphics[width=\textwidth]{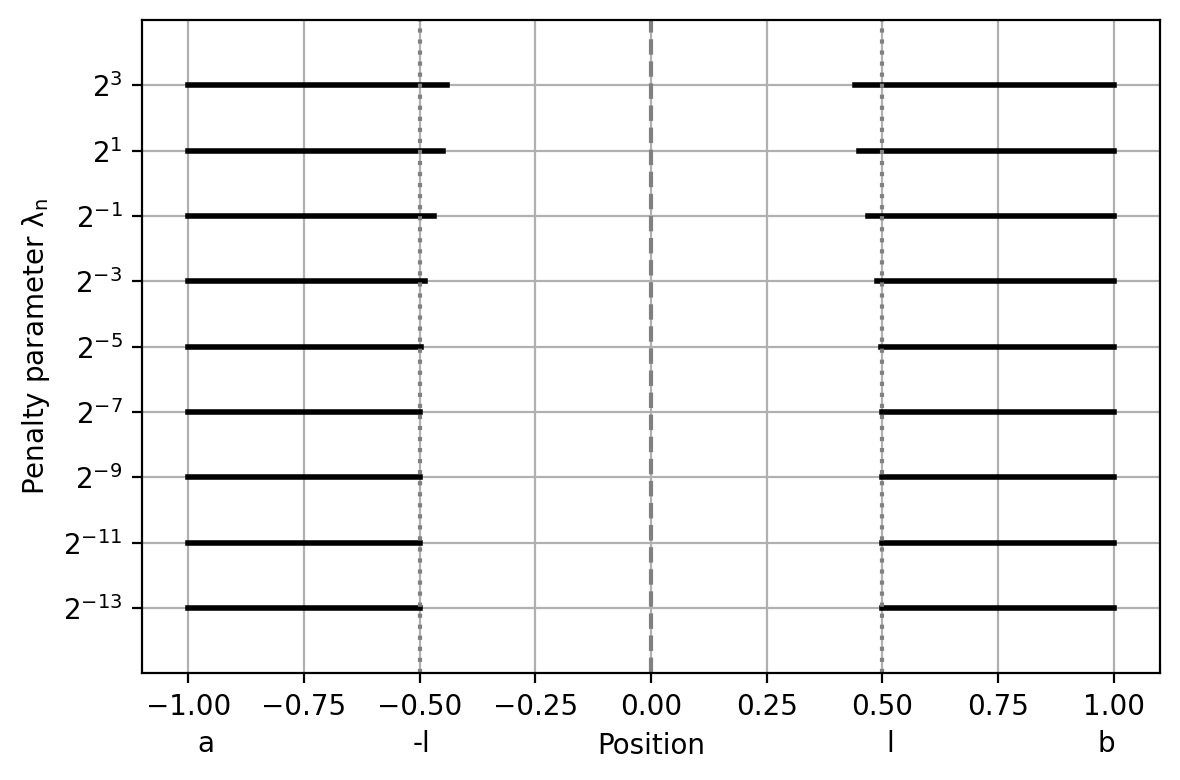}
    	    \caption{Equilibrium configuration.}
         \label{fig:convergence_displacement_1}
        \end{subfigure}
        \begin{subfigure}[b]{0.49\textwidth}
    	    \centering
    	    \includegraphics[width=\textwidth]{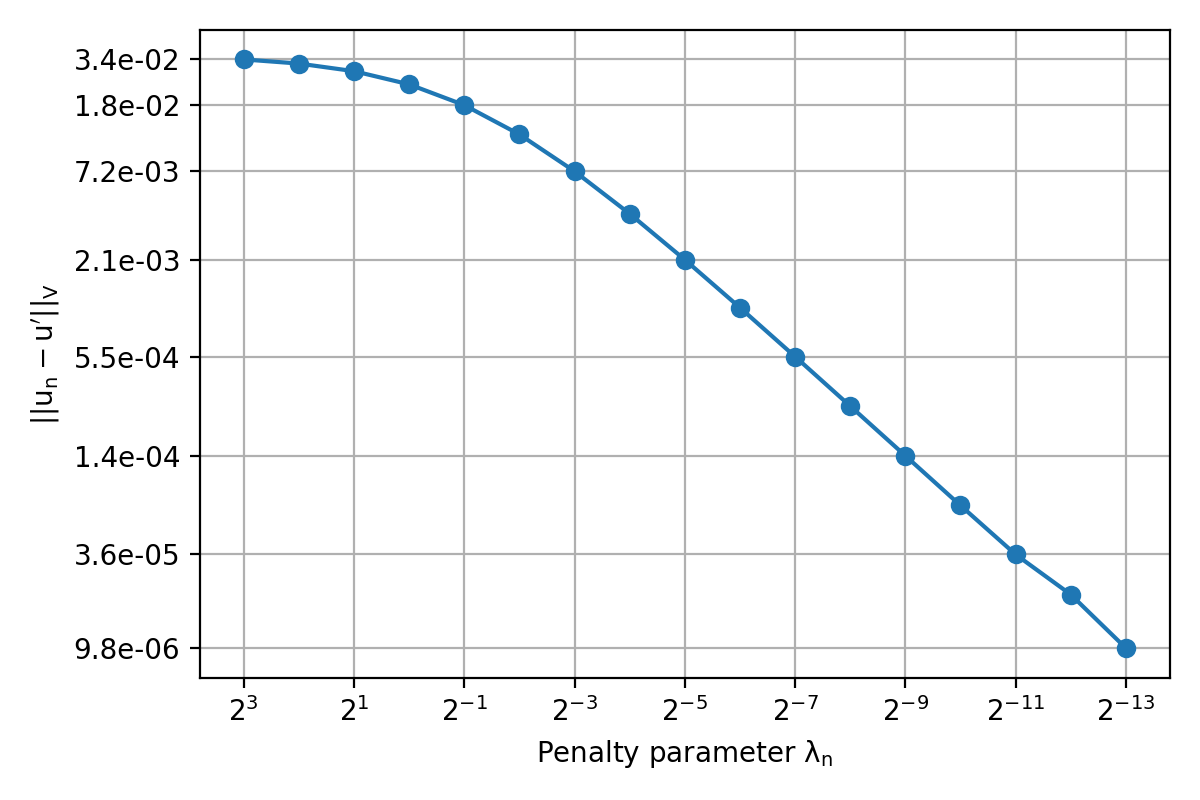}
            \caption{Norm of the approximation error.}
    	    \label{fig:convergence_error_norm_1}
        \end{subfigure}
        \begin{subfigure}[b]{0.49\textwidth}
    	    \centering
    	    \includegraphics[width=\textwidth]{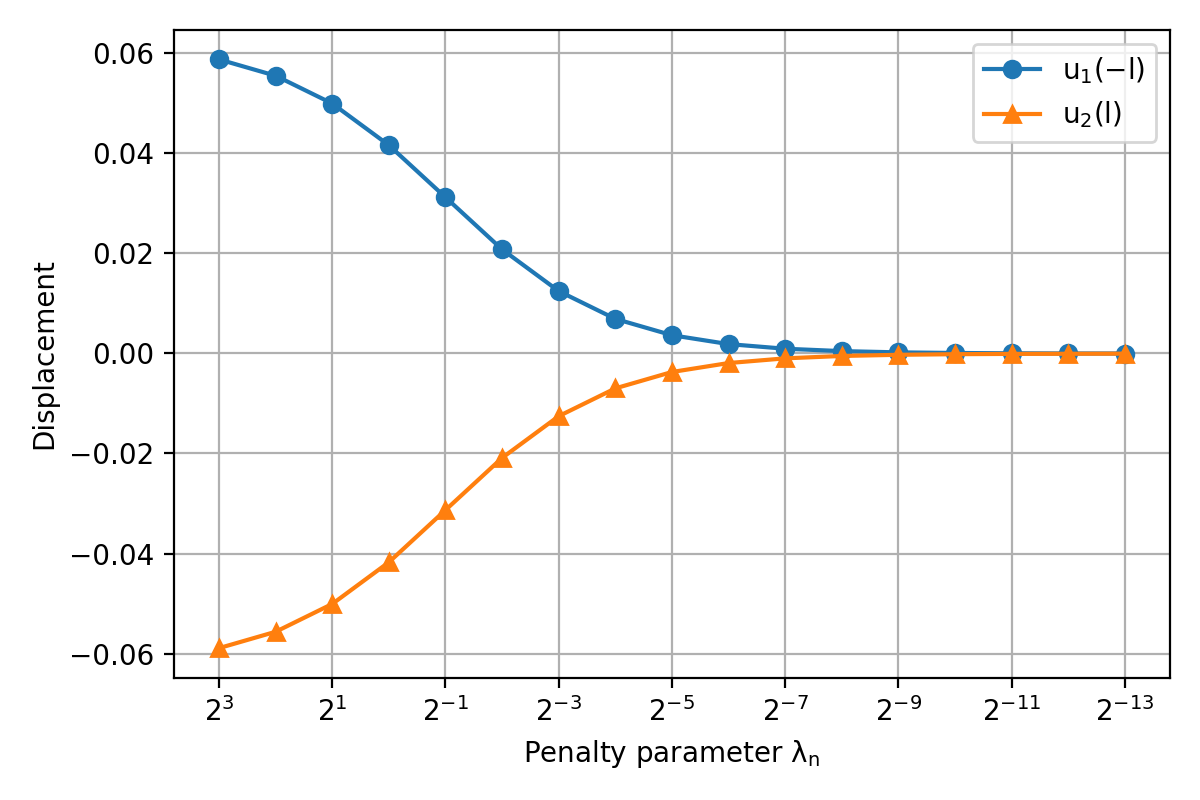}
    	    \caption{Displacement of ends of rods.}
    	    \label{fig:convergence_spring_ends_1}
        \end{subfigure}
        \begin{subfigure}[b]{0.49\textwidth}
    	    \centering
    	    \includegraphics[width=\textwidth]{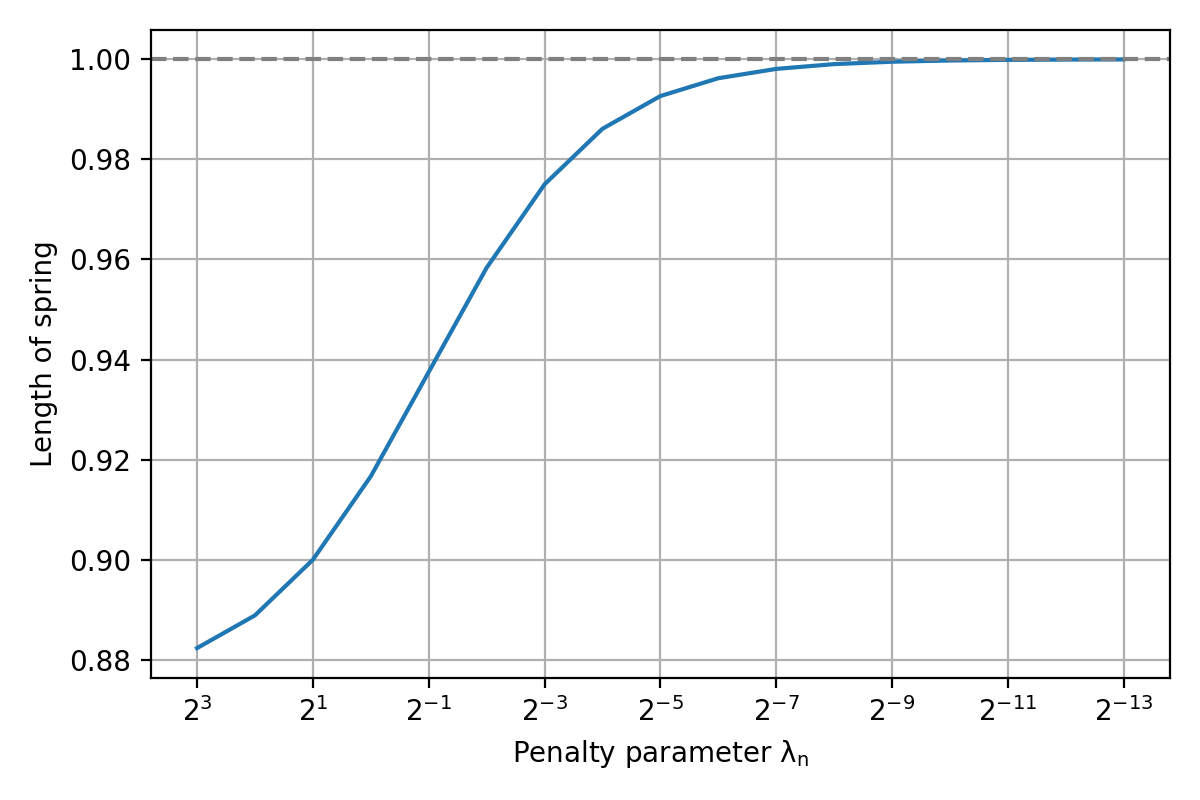}
            \caption{Length of spring.}
    	    \label{fig:convergence_spring_length_1}
        \end{subfigure}
        \medskip
        \caption{Numerical validation of the convergence $u_n\to u'$ as $\lambda_n \rightarrow 0$.}
        \label{fig:convergence_results_1}
    \end{figure}

\begin{figure}\label{fig7}
    	\vspace{-13mm}
        \centering
        \begin{subfigure}[b]{0.49\textwidth}
    	    \centering
    	    \includegraphics[width=\textwidth]{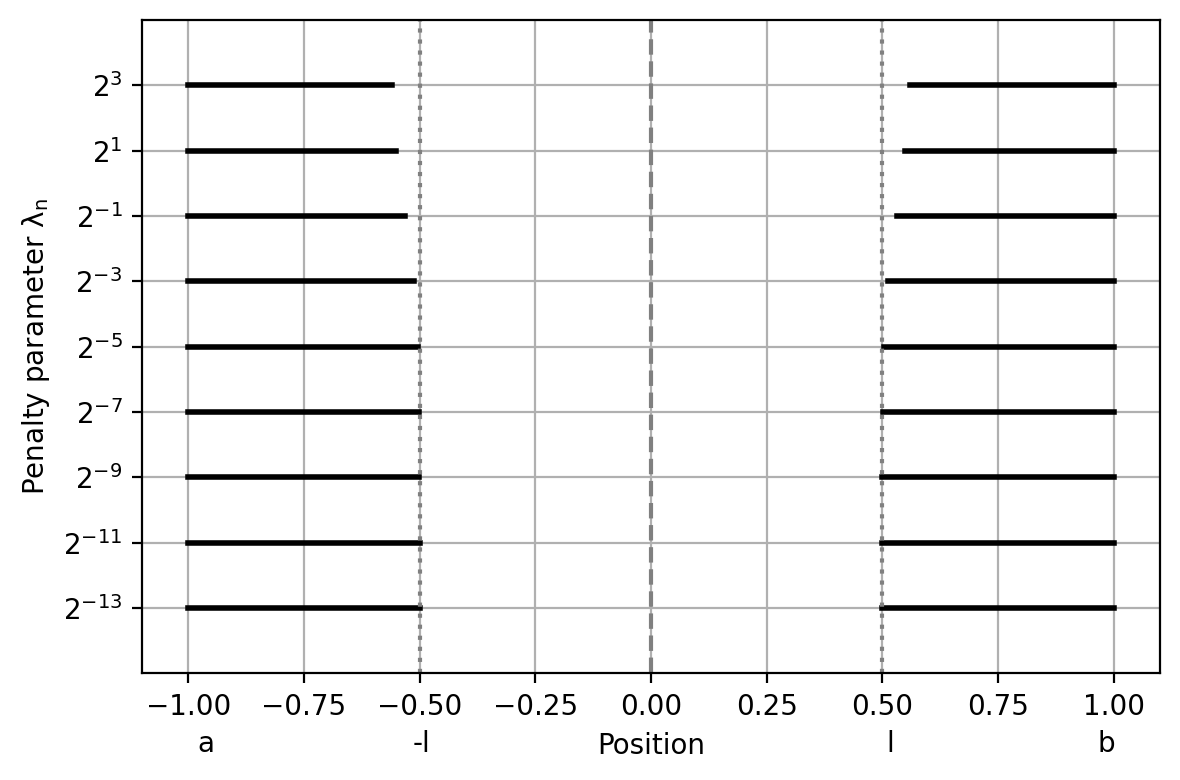}
    	    \caption{Equilibrium configuration.}
         \label{fig:convergence_displacement_2}
        \end{subfigure}
        \begin{subfigure}[b]{0.49\textwidth}
    	    \centering
    	    \includegraphics[width=\textwidth]{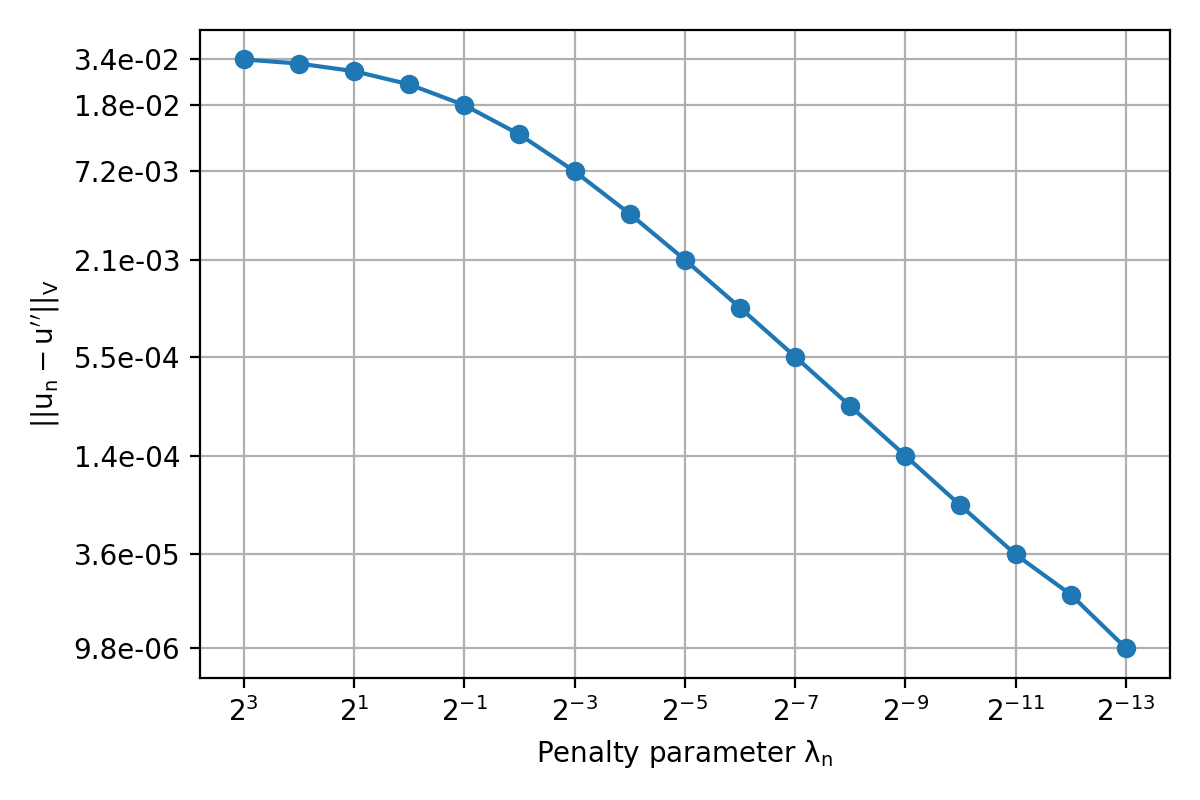}
            \caption{Norm of the approximation error.}
    	    \label{fig:convergence_error_norm_2}
        \end{subfigure}
        \begin{subfigure}[b]{0.49\textwidth}
    	    \centering
    	    \includegraphics[width=\textwidth]{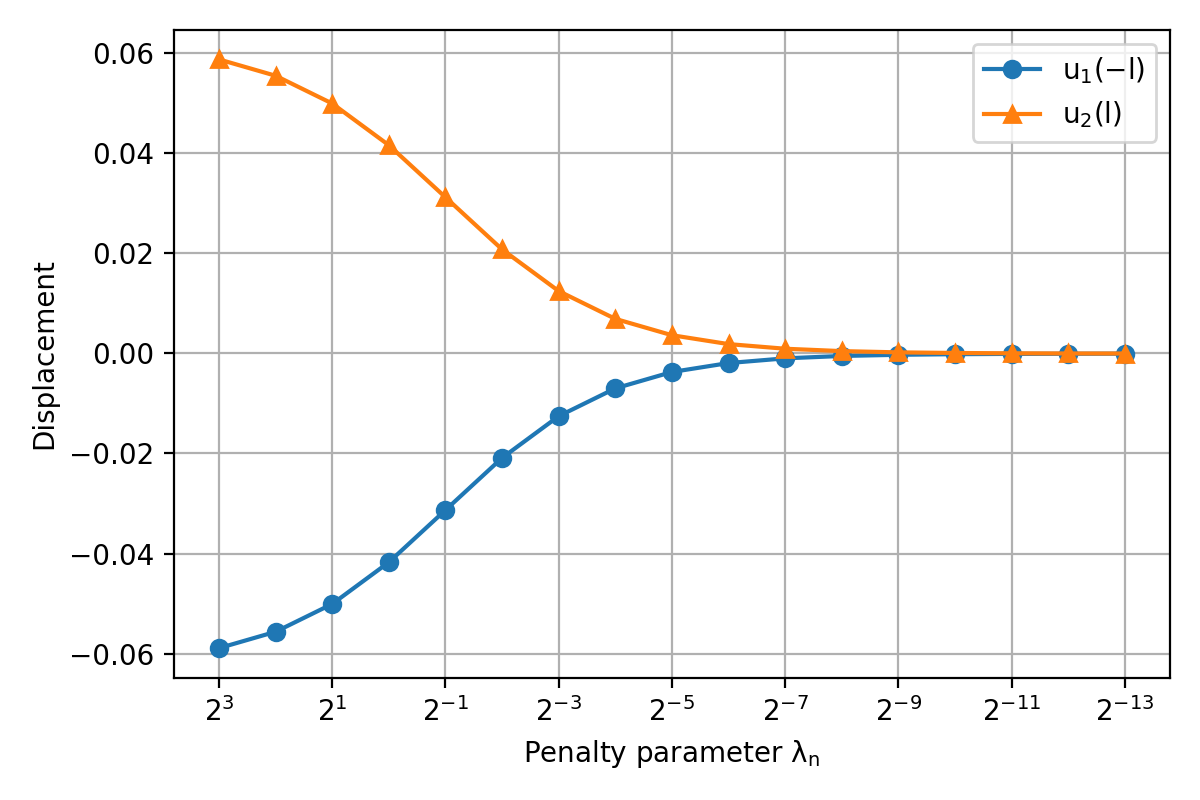}
    	    \caption{Displacement of ends of rods.}
    	    \label{fig:convergence_spring_ends_2}
        \end{subfigure}
        \begin{subfigure}[b]{0.49\textwidth}
    	    \centering
    	    \includegraphics[width=\textwidth]{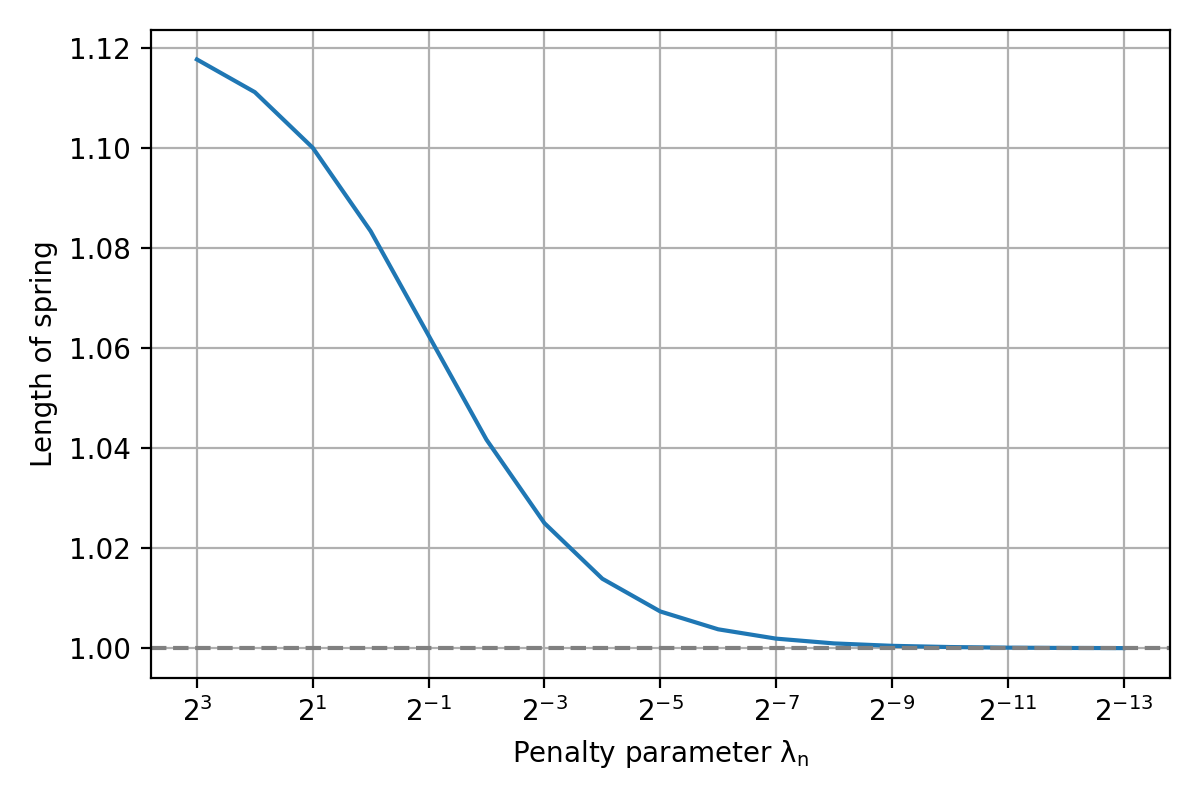}
            \caption{Length of spring.}
    	    \label{fig:convergence_spring_length_2}
        \end{subfigure}
        \medskip
        \caption{Numerical validation of the convergence $u_n\to u''$ as $\lambda_n \rightarrow 0$.}
        \label{fig:convergence_results_2}
    \end{figure}

\medskip
Our numerical results are shown in Figures \ref{fig:experiment1_figures}--\ref{fig:convergence_results_2} and are described as below.

    \begin{Experiment}\label{e1}
  {\rm  We consider the spring-rods system in which the body forces act on the rods toward the spring, i.e., $ f_1 = 1,\ \ f_2 = -1$.
    This implies the compression of the spring which pushes the ends of the rods. This situation corresponds to the one shown in Figure~\ref{fig}\,(c).
    The numerical results for this experiment show that, in the equilibrium configuration,  the length of each rods reduces as the stiffness coefficient  $k_1$ increases.
    Moreover, the magnitude of the displacements fields at the end of the rods decreases as $k_1$ increases while the magnitude of the corresponding stress vectors increases with $k_1$. The details are presented in Figure~$\ref{fig:experiment1_figures}$.}\end{Experiment}

  \vspace{-3mm}
  
  \begin{Experiment}\label{e2} 
  {\rm In the second experiment, the body forces are acting towards the rods, i.e.,
  $f_1 = -1, \ \  f_2 = 1$.    
  The spring is then extended and pulls the ends of the rods. This situation is depicted in Figure~\ref{fig}\,(b). The values of the displacements and stresses at the ends of the rods are plotted in Figure~$\ref{fig:experiment2_figures}$, for which we have similar comments to those  concerning   Figure~$\ref{fig:experiment1_figures}$.}
\end{Experiment}

\begin{Experiment}\label{e3}  
{\rm This experiment models the behavior of the system when no body force is applied to the right rod and a body force is applied to the left rod, which pushes it into the spring, i.e.,
	$f_1 = 1, \ \ f_2 = 0$.
	It can be seen that the spring transmits the force applied to the left rod, so it pushes the end of the right rod, causing its compression. 
	As the spring stiffness coefficient $k_1$ increases, the magnitude of the displacements of the ends of the rods decreases. The value of the corresponding stress at the end of the rods is increasing with $k_1$. The details are presented in 
	Figure $\ref{fig:experiment3_figures}$. }
\end{Experiment}
    
\vspace{-2mm}
\begin{Experiment}\label{e5}  
{\rm This experiment refers to the spring-rods system with equal body forces applied to the both rods, i.e.,
	$f_1 = f_2 = 1$.  
	In Figures~\ref{fig:experiment5_figures}\,(a) and  \ref{fig:experiment5_figures}\,(b) it can be seen that the displacements of the ends of the two rods are the same. Therefore, 	the length of the spring remains unchanged. This shows that
	the spring moves as if it were rigid, and it does not act on the rods with any reactive force, as shown in Figure~\ref{fig:experiment5_figures}\,(c). }
\end{Experiment}

\vspace{-2mm}
\begin{Experiment}\label{e5.5}  
	{\rm 
	The last setup of experiment is designed to show the behavior of the system in a situation
	in which the body forces in the rods, i.e. $f_1 = 6, f_2 = -6$ lead to full compression
	of the spring. At the Figure~\ref{fig:experiment5_figures}\,(a) it can be seen that for the value of spring stiffness coefficient $k_1$ below 0.5 the rods are in contact. The Figure~\ref{fig:experiment5_figures}\,(b) shows the displacements at the ends of the rods does not change for all values of $k_1$ for which rods are in contact. Furthermore, the Figure~\ref{fig:experiment5_figures}\,(c) shows the behavior of the spring stiffness, which value is lower than the value of $-p(\theta(u))$ when rods are in contact.}
	
\end{Experiment}

 \begin{Experiment}\label{e6}  
 {\rm This experiment refers to the numerical validation of the convergence result in Theorem $\ref{t2}$ a).
 	Let $u'$ denote the solution of problem $(\ref{vp})$, which models the case when the spring behaves rigidly in  compression, see $(\ref{2.4b})$.
 	We perform simulation for the solution of problem $(\ref{2m})$  with $f_1 =1$,  $f_2 = -1$, $k_1=k_2=1$, $q$ given by $(\ref{qq})$ and  $\lambda_n=\frac{1}{2^{n-3}}$,  for various values of $n\in\mathbb{N}$. The solution of this problem is denoted by $u_n$ and is obtained by considering an approximating method for the solution of a minimization Problem~$\hat{\cP}_V^n$ similar to Problem~$\hat{\cP}_V$. We see that for $n$ large enough the length of the spring approaches the initial length $2l=1$, as shown in Figure $\ref{fig:convergence_results_1}$\,(a), (d).  In Figure $\ref{fig:convergence_results_1}$\,(b)  we plot the values of the norm $\|u_n-u'\|_V$ and note that it converges to zero as $n\to\infty$.
 	The evolution of displacements at the ends of the rods are  plotted  in Figure $\ref{fig:convergence_results_1}$\,(c). It follows from there that for $n$ large enough the corresponding displacements are zero and, in exchange, the length of the spring is $2l$. All these results show that for $\lambda_n$ small enough, the solution of the penalized Problem~$\cP_V^n$ approaches the solution of Problem~$\cP_V$, in which the spring has a rigid behavior in compression. They represent a numerical validation of the convergence result in Theorem $\ref{t2}$ a).}
 \end{Experiment}

    \filbreak

    \begin{Experiment}\label{e7} 
    {\rm   The experiment refers to the numerical validation of the convergence result in Theorem $\ref{t2}$ b). 
    	Let $u''$ denote the solution of problem $(\ref{vpp})$, which models the case when the spring behaves rigidly in  extension, see $(\ref{2.4c})$.
    	We perform simulation for the solution of problem $(\ref{1m})$  with $f_1 =1$,  $f_2 = -1$, $k_1=k_2=1$, $q$ given by \eqref{qqn} and  $\lambda_n=\frac{1}{2^{n-3}}$,  for various values of $n\in\mathbb{N}$. The solution of this problem is denoted by $u_n$
    	and, again, it is obtained by considering an approximating method for the solution of a minimization Problem~$\hat{\cP}_V^n$ similar to Problem~$\hat{\cP}_V$. Our results are presented in Figure~\ref{fig:convergence_results_2} and have a similar interpretation as those of the previous experiment:
    	for $n$ large enough the length of the spring approaches the initial length $2l=1$, the norm $\|u_n-u''\|_V$ converges to zero as $n\to\infty$, the evolution of displacements at the ends of the rods vanish for n large enough.  All these results show that for  $\lambda_n$ small enough the solution of the penalized Problem~$\cP_V^n$ approaches the solution of Problem~$\cP_V$, in which the spring has a rigid behavior in extension. They represent a numerical validation of the convergence result in Theorem~\ref{t2} b).}
    \end{Experiment}

\vskip 2mm
\noindent {\bf Acknowledgments}\\
The project has received funding from the European Union's Horizon 2020 Research and Innovation Programme under the Marie Sklo\-do\-wska-Curie grant agreement no.\ 823731 CONMECH.
The second author is supported by the project financed by the Ministry of Science and Higher Education of Republic of Poland under Grant No. 440328/PnH2/2019, and in part from National Science Center, Poland, under project OPUS no. 2021/41/B/ST1/01636.

\end{document}